\documentclass[11 pt]{article}
\usepackage{verbatim}
\usepackage{amsmath}
\usepackage{latexsym}
\usepackage{amssymb}
\usepackage{amscd}
\usepackage{amsthm}
\newtheorem{Thm}{Theorem}[section]

\newtheorem{Lem}[Thm]{Lemma}

\newtheorem{Prob}{Problem}

 \makeatletter
\newfont{\footsc}{cmcsc10 at 8truept}
\newfont{\footbf}{cmbx10 at 8truept}
\newfont{\footrm}{cmr10 at 10truept}
\title{Number of Irreducible Polynomials and Pairs of Relatively Prime Polynomials in Several Variables over Finite Fields}

\author
{Xiang-dong Hou\footnote{Department of Mathematics,  University of South Florida, 
Tampa, FL 33620, USA; Email: xhou@cas.usf.edu} 
\ and Gary L. Mullen\footnote{Department of
Mathematics, The Pennsylvania State University, University Park,
PA 16802, USA; Email: mullen@math.psu.edu} }

\date{}

\begin{document}
\maketitle

\newcommand{\bc}[2]{{{#1}\choose{#2}}}
\def\mod{\mathop{\mathrm{mod}}}
\def\deg{\mathop{\mathrm{deg}}}

\abstract{We discuss several enumerative results for irreducible polynomials of a given degree and pairs 
of relatively prime polynomials of given degrees in 
several variables over finite fields.
Two notions of degree, the {\em total degree} and the {\em vector degree}, 
are considered. We show that the number of irreducibles can be 
computed recursively by degree and that the number of relatively prime 
pairs can be expressed in terms of the number of irreducibles. We 
also obtain asymptotic formulas for the number of irreducibles and 
the number of relatively prime pairs. The asymptotic formulas for the 
number of irreducibles generalize and improve several previous results by 
Carlitz, Cohen and Bodin. 
}

\section{Introduction} 

Let $\Bbb F_q$ be the finite field with $q$ elements. In this paper we consider two problems:

\begin{Prob}~\label{Pb1.1}
Count the number of irreducible polynomials of a given degree in $\Bbb F_q[x_1,\dots,x_k]$.
\end{Prob}

\begin{Prob}~\label{Pb1.2}
Count the number of pairs of relatively prime polynomials of given 
degrees in $\Bbb F_q[x_1,\dots,x_k]$.
\end{Prob}

When $k=1$, both problems have been solved. The following formula for the number $I(m)$  of monic irreducible 
polynomials of degree $m$ in $\Bbb F_q[x]$ is well known (see \cite{LidNie}):
\begin{equation}\label{Nqm}
I(m) = \frac {1} {m} \sum _{d \vert m} \mu (d) q^{m/d}.
\end{equation} 

In \cite{Cor} the authors show that the number 
of pairs of polynomials $f(x)$ and $g(x)$ of degree $m$ 
over the binary field $\Bbb F_2$  with greatest common
divisor $(f,g)=1 $ is the same as the number of  pairs of polynomials of degree $m$ in which
$(f,g)\ne 1$. The authors also asked for a \lq\lq nice simple bijection that proves this
result."
In \cite {Rei}  a bijection using \lq\lq resultant matrices" is found. 
More recently in \cite{BenBen},
using the Euclidean Algorithm, 
the authors exhibit a more natural bijection between  pairs of
binary polynomials $f(x)$ and $g(x)$ of degree $m$  with the greatest common
divisor $(f,g)=1 $ and pairs of polynomials of degree $m$ with
$(f,g)\ne 1$. In fact, the following result of \cite{BenBen} answers more than the $k=1$ case of Problem~\ref{Pb1.2}.

\begin{Thm} {\rm(\cite[Corollary 5]{BenBen})}\label{T1.3}
Let $(0,\dots,0)\ne (d_1,\dots,d_s)\in\Bbb N^s$ and let 
$f_i(x)\in \Bbb F_q[x]$ be a randomly chosen polynomial of degree $d_i$. Then the probability that
$\text{\rm gcd}(f_1,\dots,f_s)=1$ is $1-\frac 1{q^{s-1}}$.
\end{Thm}
 
When $k\ge 2$, the situation for both problems is quite different. First of all, there are no known closed formulas for the numbers in the two problems.
Formula~\eqref{Nqm} depends on the fact that the polynomial $x^{q^m} -x$ is the product of all monic irreducible polynomials over $\Bbb F_q$ of degree
$d$  where $d$ divides $m$; see Theorem 3.20 of \cite{LidNie}. Unfortunately, 
there is no known analogous result for polynomials in two or more variables.

Before we proceed, it should be pointed out that in $\Bbb F_q[x_1,\dots,x_k]$ with $k\ge 2$, there are two notions of degree. Let $0\ne f(x_1,\dots,x_k)\in
\Bbb F_q[x_1,\dots,x_k]$. The {\em total degree} of $f$, denoted by $\deg f$, is the degree of the polynomial $f(tx_1,\dots,tx_k)$ in $t$ over $\Bbb F_q[x_1,\dots,x_k]$. The {\em vector degree} of $f$, denoted by $\text{Deg}\,f$, is the $k$-tuple $(\deg_{x_1}f,\dots,\deg_{x_k}f)$. Thus each of the above two problems has a total degree version and a vector degree version; treatments of the two versions are not entirely the same.

Carlitz \cite{Car63} studied Problem~\ref{Pb1.1} with total 
degree and obtained an asymptotic formula for the number of 
irreducible $f\in\Bbb F_q[x_1,\dots,x_k]$ with $\deg f=m$ as $m\to \infty$. 
Recently, Bodin \cite{Bod} improved Carlitz's result by providing the next 
term in Carlitz's asymptotic formula.
Paper \cite{Bod} also gives a recursive formula for computing the number of irreducible $f\in\Bbb F_q[x_1,\dots,x_k]$ with $\deg f=m$.
The study of Problem~\ref{Pb1.1} with vector degree started with Carlitz
 \cite{Car65} in which he obtained asymptotic results for the number of 
 irreducible $f\in\Bbb F_q[x_1,x_2]$ with $\text{Deg}\,f=(m_1,m_2)$. This was later generalized to an arbitrary number ($\ge 2$) of variables by Cohen \cite{Coh}; for Cohen's further work on the topic, see \cite{Coh2,Coh3}.

There is a fundamental difference between irreducible polynomials in one variable and those in several variables. When $k=1$, as $m\to \infty$, almost all polynomials 
of degree $m$ are
{\em reducible}; this follows easily from \eqref{Nqm}. However, when $k\ge 2$, as $m\to\infty$, almost all polynomials of total degree $m$ are {\em irreducible}; see 
\cite[Theorem 7]{Bod} or Theorem~\ref{T2.2} of the present 
paper. When $k\ge 2$ and $m_1,\dots,m_{k-2}$ are fixed, as 
$m_{k-1},m_k\to\infty$, almost all polynomials of degree $(m_1,\dots,m_k)$ are {\em irreducible}; see 
\cite[Theorem~1]{Coh} and Theorem~\ref{T5.3} of the present paper.

Problem~\ref{Pb1.2} with $k\ge 2$ was the initial motivation for our work. We want to see to what extent an analogous kind of result on the number relatively prime pairs might hold in several variables.
The only published result we are aware of is Corollary 12 of \cite{Coh3}. It states that the proportion of the relative prime pairs of polynomials of degree $(m_1,\dots,m_k)$ tends to $1-q^{1-2(m_1+1)\cdots(m_{k-1}+1)}$ as $m_k\to\infty$ (with $m_1,\dots, m_{k-1}$ fixed). We will consider pairs of polynomials in several variables not necessarily of the same degree.  
We find  that unlike the one variable case, 
almost all pairs of polynomials in several variables are relatively prime.

We now introduce some basic notation.  
Let $\mathcal N_k=\Bbb F_q[x_1,\dots,x_k]/\!\sim$, where 
$f\sim g$  if $f=cg$ for some $c\in\Bbb F_q^\times$. Elements in $\mathcal N_k$ are normalized polynomials in $k$ variables
which correspond to monic polynomials in one variable. 
For $m,n\in\Bbb N$ and $\frak m=(m_1,\dots,m_k)$, $\frak n=(n_1,\dots,n_k)\in\Bbb N^k$, let
\begin{gather*}
\mathcal N_k(m)=\{f\in \mathcal N_k:\deg f=m\},\qquad N_k(m)=|\mathcal N_k(m)|,\\
\mathcal N_k(\frak m)=\{f\in \mathcal N_k:\text{Deg}\, f=\frak m\},\qquad N_k(\frak m)=|\mathcal N_k(\frak m)|,\\
\mathcal I_k(m)=|\{f\in \mathcal N_k(m): \text{$f$ is irreducible}\},\qquad I_k(m)=|\mathcal I_k(m)|,\\
\mathcal I_k(\frak m)=|\{f\in \mathcal N_k(\frak m): \text{$f$ is irreducible}\},\qquad I_k(\frak m)=|\mathcal I_k(\frak m)|,\\
P_k(m;n)=|\{(f,g)\in \mathcal N_k(m)\times \mathcal N_k(n):\text{gcd}(f,g)=1\}|,\\
P_k(\frak m;\frak n)=|\{(f,g)\in \mathcal N_k(\frak m)\times \mathcal N_k(\frak n):\text{gcd}(f,g)=1\}|.
\end{gather*}
   
We next summarize the contributions of the current paper. A recursive 
formula for $I_k(m)$ has been given in \cite{Bod}. We show that a similar 
formula holds for $I_k(\frak m)$. (In fact, the recursive formula works 
for any grading of $\Bbb F_q[x_1,\dots,x_k]$ by a partially ordered monoid; see \cite[Ch. II, \S 11.2]{Bou89}. For example, one can grade 
$\Bbb F_q[x_1,\dots,x_k]$ by total degrees on several subsets of $\{x_1,\dots,x_k\}$.) We provide formulas 
for $P_k(m;n)$ in terms of $I_k(i)$ ($i\le\min\{m,n\}$) and for 
$P_k(\frak m;\frak n)$ in terms of $I_k(\frak i)$ 
($\frak i\le \frak m,\frak n$). We obtain asymptotic formulas for $I_k(m)$, 
$I_k(\frak m)$, $P_k(m;n)$ and $P_k(\frak m;\frak n)$.
The asymptotic formula for $I_k(m)$ (as $m\to \infty$) is an 
expansion of $I_k(m)$ with explicit terms and accurate up to $O(q^{\binom{m-t-1+k}k})$ for any $t\ge 0$. 
The results of \cite{Car63} and \cite{Bod} are special cases of this expansion with one term and two terms respectively. 
Our asymptotic formula for $I_k(\frak m)$ is an improvement of the one in \cite{Coh}.

The paper is organized as follows. Section 2 deals with the total degree version of Problem~\ref{Pb1.1}. In section 3 we briefly describe an algorithm for computing the gcd of two polynomials in $\Bbb F_q[x_1,\dots,x_k]$. Section 4 is devoted to the total degree version of Problem~\ref{Pb1.2}. The vector degree version of both Problems~\ref{Pb1.1} and \ref{Pb1.2} is discussed in section 5. Appendix A contains the deferred proof of Lemma~\ref{L0.3} which is rather lengthy and technical. Appendix B contains several tables of values of the functions $I_k(m)$, $I_k(\frak m)$, $P_k(m;n)$, $P_k(\frak m;\frak n)$. 

We conclude this section with a quick review of the M\"obius inversion formula which is a basic tool of this paper. We refer the reader to \cite{Ben75}
for more details on the subject. Let $(X,\le)$ be a partially ordered set such that for all $x,y\in X$, the interval $[x,y]=\{z\in X:x\le z\le y\}$ is finite. The M\"obius function of $(X,\le)$ is the function $\mu:X\times X\to \Bbb Z$ such that
\[
\sum_{z\in[x,y]}\mu(x,z)=\begin{cases}
1&\text{if}\ x=y,\cr
0&\text{if}\ x\ne y.
\end{cases}
\]
Let $A$ be an abelian group and let $N_{=}:X\to A$ be a function. Fix $l,m\in X$ and for $x\in X$ define
\begin{gather*}
N_{\ge}(x)=\sum_{y\in [x,m]}N_{=}(y),\\
N_{\le}(x)=\sum_{y\in [l,x]}N_{=}(y).
\end{gather*}
Then we have
\[
N_{=}(x)=\sum_{y\in [x,m]}\mu(x,y)N_{\ge}(y)\qquad\text{for all $x\in X$ with $x\le m$}
\]
and 
\[
N_{=}(x)=\sum_{y\in [l,x]}\mu(y,x)N_{\le}(y)\qquad\text{for all $x\in X$ with $x\ge l$}.
\]
If $(X,\le)$ has a minimum element $1$, $\mu(1,x)$ is usually denoted by $\mu(x)$.


\section{Number of irreducible polynomials in several variables} 


\subsection{Recursive formula for $I_k(m)$}

A recursive formula for $I_k(m)$ can be found in \cite{Bod}. The only new contribution in this subsection is some computational and numerical results.

Let $k\ge 1$ and $m\ge 0$. It is easy to see that 
\[
N_k(m)=q^{\binom{m+k-1}k}\frac{q^{\binom{m+k-1}{k-1}}-1}{q-1}=\frac{q^{\binom{m+k}k}-q^{\binom{m+k-1}k}}{q-1}.
\]
Unique factorization in $\Bbb F_q[x_1,\dots,x_k]$ implies that
\[
\sum_{1a_1+2a_2+\cdots+ma_m=m}\binom{I_k(1)+a_1-1}{a_1}\cdots\binom{I_k(m)+a_m-1}{a_m}=N_k(m).
\]
(In the above sum, $\binom{I_k(i)+a_i-1}{a_i}$ is the number of products of $a_i$ (not necessarily distinct) elements from $\mathcal I_k(i)$.)
This allows us to compute $I_k(m)$ recursively. Starting with $I_k(0)=0$, we have for $m>0$
\begin{equation}\label{2a}
I_k(m)=N_k(m)-\sum_{1a_1+2a_2+\cdots+(m-1)a_{m-1}=m}\binom{I_k(1)+a_1-1}{a_1}\cdots\binom{I_k(m-1)+a_{m-1}-1}{a_{m-1}}.
\end{equation}

We next provide explicit formulas for $I_k(m)$ with $m\le 3$ and $I_2(m)$ with $m\le 10$, obtained from \eqref{2a} using Mathematica~\cite{Wol}.
\[
I_k(0)=0.
\]
\[
I_k(1)=\frac q{q-1}(q^k-1).
\]
\[
I_k(2)=\frac q{2(q-1)^2}\bigl[ 2(q-1)q^{\frac 12 k(k+3)}-q^{2k+1}-q^{k+1}+3q^k-1\bigr].
\]
\[
\begin{split}
I_k(3)=\,&\frac q{3(q-1)^3}\bigl[3(q-1)^2q^{\frac 16k(k^2+6k+11)}-3(q-1)q^{\frac 12 k(k+3)}(q^{k+1}-1)\cr
&+q^{3k+2}-3q^{2k+1}-q^{k+2}+5q^{k+1}-q^k-2q+1\bigr].
\end{split}
\]
\[
I_2(0)=0.
\]
\[
I_2(1)=q (q + 1).
\]
\[
I_2(2)=\frac{1}{2} (q - 1) q (2 q + 1) \left(q^2 + q + 1\right).
\]
\[
I_2(3)=\frac{1}{3} (q - 1) q (q + 1) \left(3 q^6 + 3 q^5 + 3 q^4 + q^3 + q + 1\right).
\]
\[
I_2(4)=\frac{1}{4} (q - 1) q^2 \left(
      q^2 + q + 1\right) \left(4 q^9 + 4 q^8 + 4 q^7 + 
        4 q^6 - 2 q^5 - 4 q^4 + q^3 + 2 q^2 - 1\right).
\]
\[
I_2(5)=\frac{1}{5} (q - 1) q (q + 1) \left(q^2 + q + 1\right) \left(
    5 q^{15} + 5 q^{13} + 5 q^{12} - 5 q^9 - 5 q^8 + 5 
    q^6 + q^5 - q^4 + q^3 + 1\right).
\]
\[
\begin{split}
I_2(6)=\,&\frac{1}{6} (q - 1) q (q + 1) \left(q^2 + q + 1\right) \Bigl(
    6 q^{22} + 6 q^{20} + 6 q^{
        19} + 6 q^{18} + 6 q^{16} - 6 q^{15} - 6 q^{14}\cr
& - 9 q^{
      13} - 9 q^{12} + 11 q^{10} + 7 q^9 + 2 q^8 - 5 q^7 - q^6 + q^5 - 
      q^4 - 2 q^3 + 2 q^2 - q - 1\Bigr).
\end{split}
\]
\[
\begin{split}
I_2(7)=\,&\frac{1}{7} (q - 1) q (q + 1) \left(q^2 + q + 1\right) \Bigl(
    7 q^{30} + 7 q^{28} + 7 q^{
        27} + 7 q^{26} + 7 q^{25} + 7 q^{24} - 14 q^{20}\cr
& - 7 q^{19} - 14 q^{18} - 14 q^{17} - 7 q^{16} + 7 q^{15} + 21 q^{14} + 
      21 q^{13} + 7 q^{12} - 14 q^{11} - 14 q^{10}\cr
& + q^9 + 6 q^8 + q^7 + 
      q^3 + 1\Bigr).
\end{split}
\]
\[
\begin{split}
I_2(8)=\,&\frac{1}{8} (q - 1) q^4 (q + 1) \left(q^2 + q + 1\right) \Bigl(8 
    q^{36} + 8 q^{34} + 8 q^{33} + 
        8 q^{32} + 8 q^{31} + 16 q^{30} + 8 q^{28}\cr
& - 8 q^{25} - 8 q^{24} - 16 
      q^{23} - 8 q^{22} - 16 q^{21} - 20 q^{20} - 12 q^{19} - 8 q^{18} + 4 
      q^{17} + 24 q^{16}\cr
& + 44 q^{15} + 40 q^{14} + 4 q^{13} - 30 q^{12} - 
      46 q^{11} - 14 q^{10} + 20 q^9 + 21 q^8 - q^7 - 9 q^6 \cr
&- 2 q^5 - 2 
      q^4 + 3 q^2 + q - 1\Bigr).
\end{split}
\]
\[
\begin{split}
I_2(9)=\,&\frac{1}{9} (q - 1) q^3 (q + 1) \left(q^2 + q + 1\right) \Bigl(9 
    q^{47} + 9 q^{45} + 9 q^{44} + 
        9 q^{43} + 9 q^{42} + 18 q^{41} + 9 q^{40}\cr
& + 9 q^{39} + 9 q^{38} - 9 
      q^{34} - 18 q^{33} - 9 q^{32} - 18 q^{31} - 9 q^{30} - 18 q^{29} - 18 
      q^{28} - 18 q^{27}\cr
& - 9 q^{26} - 18 q^{25} - 9 q^{24} + 18 q^{23} + 
      36 q^{22} + 63 q^{21} + 66 q^{20} + 51 q^{19} - 15 q^{18} - 87 q^{
      17}\cr
& - 102 q^{16} - 42 q^{15} + 48 q^{14} + 78 q^{13} + 30 q^{12} - 26 
      q^{11} - 28 q^{10} + q^9 + 9 q^8 + q^5 - q^4\cr
& + q^3 - 2 q^2 - q + 1\Bigr).
\end{split}
\]
\[
\begin{split}
I_2(10)=\,&\frac{1}{10} (q - 1) q (q + 1) \left(q^2 + q + 1\right) \Bigl(10 q^{60} + 10 
          q^{58} + 10 q^{57} + 10 q^{56} + 10 
    q^{55} + 20 q^{54}\cr
& + 10 q^{53} + 20 q^{52} + 10 q^{51} + 
      10 q^{50} + 10 q^{48} - 10 q^{47} - 10 q^{45} - 20 q^{44} - 20 q^{43}\cr
& - 10 q^{42} - 20 q^{41} - 10 q^{40} - 20 q^{39} - 20 q^{38} - 10 q^{37} - 
      10 q^{36} - 15 q^{35} - 25 q^{34}\cr
& - 20 q^{33} - 5 q^{32} + 20 q^{31} + 
      25 q^{30} + 65 q^{29} + 100 q^{28} + 105 q^{27} + 
      75 q^{26} - 5 q^{25}\cr
& - 90 q^{24} - 195 q^{23} - 
      195 q^{22} - 55 q^{21} + 132 q^{20} + 213 q^{19} + 117 q^{18} - 60 q^{
      17}\cr
& - 125 q^{16} - 54 q^{15} + 36 q^{14} + 29 q^{13} - 3 q^{12} - 10 q^{
      11} + 6 q^9 + 6 q^8 - q^7 - 3 q^6\cr
& - 2 q^5 + 2 q^4 - q^3 + q^2 - q - 1\Bigr).
\end{split}
\]

We observe that $I_2(m)$ is a monic polynomial of degree $m(m+3)/2$ in $q$ and more generally
$I_k(m)$ is a monic polynomial of degree $\binom{m+k}k-1$ in $q$. This fact can be easily proved by induction. 
A table of values of $I_2(m)$ with $m\le 10$ and $q=2,3,4,5$ is given in the appendix. The values of $I_2(m)$ with $q=2$ and $m\le 10$ have been given in \cite{Bod}.


\subsection{Asymptotic formula for $I_k(m)$}

\begin{Lem}\label{L2.1}
Let $k\ge 2$ and $t\ge 0$ be fixed integers and also let $q$ be fixed. Then
\[
\sum_{1a_1+2a_2+\cdots+(m-t)a_{m-t}=m}\binom{I_k(1)+a_1-1}{a_1}\cdots\binom{I_k(m-t)+a_{m-t}-1}{a_{m-t}}=
O(q^{\binom{m-t+k}k}),
\]
where the $O$ concerns only the variable $m$ and
the constant in $O(q^{\binom{m-t+k}k})$ depends only on $q$, $k$ and $t$.
\end{Lem}

\begin{proof}
Assume $m>3t$. 
Let
\[
\mathcal F=\{f\in\mathcal N_k(m): \text{all irreducible factors of $f$ have deg}\ \le m-t\}.
\]
The sum in Lemma~\ref{L2.1} is $|\mathcal F|$. We claim that every $f\in\mathcal F$ can be written as $f=f_1f_2$ with $t<\deg f_1\le m-t$. Suppose to the contrary that $f\in\mathcal F$ does not allow such a factorization. Then all irreducible factors of $f$ have $\deg\le t$. Hence $f$ has a factor $f_1$ with $\frac m2-\frac t2\le \deg f_1\le \frac m2+\frac t2$. So we must have $\frac m2-\frac t2\le t$ or $\frac m2+\frac t2>m-t$, both of which are false since $m>3t$. So the claim is proved. Now we have
\[
|\mathcal F|\le\sum_{t<d\le m-t}N_k(d)N_k(m-d)\le 2\sum_{t\le d\le \frac m2}N_k(d)N_k(m-d)\le 2\sum_{t\le d\le \frac m2}q^{\binom{d+k}k+\binom{m-d+k}k}.
\]
Let $c_d=\binom{d+k}k+\binom{m-d+k}k$. Then 
for $t<d\le \frac m2$,
\[
\begin{split}
c_d-c_{d-1}=\,&\binom{d+k}k+\binom{m-d+k}k-\binom{d-1+k}k-\binom{m-(d-1)+k}k\cr
=\,&\binom{d-1+k}{k-1}-\binom{m-d+k}{k-1}\le -1\qquad\text{(since $d-1+k<m-d+k$)}.
\end{split}
\]
Thus by induction, $c_d-c_t\le -(d-t)$. So we have
\[
\begin{split}
|\mathcal F|\,&\le 2\sum_{t\le d\le \frac m2}q^{c_d}\le 2\sum_{t\le d\le \frac m2}q^{c_t-(d-t)}
= 2\sum_{t\le d\le \frac m2}q^{\binom{t+k}k+\binom{m-t+k}k-(d-t)}\cr
&=q^{\binom{m-t+k}k}\cdot 2 q^{\binom{t+k}k}\sum_{t\le d\le \frac m2}q^{-(d-t)}=O(q^{\binom{m-t+k}k}).
\end{split}
\]
\end{proof}

\begin{Thm}\label{T2.2}
Let $k\ge 2$ and $t\ge 0$ be fixed integers and also let $q$ be fixed. Then as $m\to\infty$,
\begin{equation}\label{9A}
I_k(m)=\sum_{i=0}^t\alpha_i N_k(m-i)+O(q^{\binom{m-t-1+k}k}),
\end{equation}
where the $O$ concerns only the variable $m$ and the sequence $\alpha_i$ is given by
\[
\begin{cases}
\alpha_0=1,\cr
\alpha_i=-N_k(i)\alpha_0-\cdots-N_k(1)\alpha_{i-1},&i>0.
\end{cases}
\]
\end{Thm}

\noindent{\bf Note.}
The recursive formula for $\alpha_i$ in Theorem~\ref{T2.2} is equivalent to
\begin{equation}\label{10A}
\sum_{i=0}^\infty\alpha_i x^i=\Bigl(\sum_{i=0}^\infty N_k(i)x^i\Bigr)^{-1}.
\end{equation}
From \eqref{10A} one can derive the following explicit formula for $\alpha_i$:
\[
\alpha_i=\sum_{1a_1+\cdots+ia_i=i}\frac{(a_1+\cdots+a_i)!}{a_1!\cdots a_i!}(-1)^{a_1+\cdots+a_i}N_k(1)^{a_1}\cdots N_k(i)^{a_i},\qquad i>0.
\]

\begin{proof}[Proof of Theorem~\ref{T2.2}] Use induction on $t$. By \eqref{2a} and Lemma~\ref{L2.1} we have
\[
I_k(m)=N_k(m)+O(q^{\binom{m-1+k}k}).
\]
So the conclusion holds for $t=0$. Now assume $t>0$. When $m$ is large, \eqref{2a} and Lemma~\ref{L2.1} give
\[
\begin{split}
I_k(m)=\,&N_k(m)-N_k(1)I_k(m-1)-\cdots-N_k(t)I_k(m-t)\cr
&-\sum_{1a_1+2a_2+\cdots+(m-t-1)a_{m-t-1}=m}\binom{I_k(1)+a_1-1}{a_1}\cdots\binom{I_k(m-t-1)+a_{m-t-1}-1}{a_{m-t-1}}\cr
=\,&N_k(m)-N_k(1)I_k(m-1)-\cdots-N_k(t)I_k(m-t)+O(q^{\binom{m-t-1+k}k}).
\end{split}
\]
By the induction hypothesis, the last expression equals
\[
\begin{split}
&N_k(m)-N_k(1)\bigl(\alpha_0 N_k(m-1)+\cdots+\alpha_{t-1} N_k(m-t)\bigr)-\cdots-N_k(t)\alpha_0 N_k(m-t)\cr
&+O(q^{\binom{m-t-1+k}k})\cr
=\,&N_k(m)-N_k(1)\alpha_0 N_k(m-1)-\cdots-\bigl(N_k(t)\alpha_0+\cdots+N_k(1)\alpha_{t-1}\bigr)N_k(m-t)\cr
&+O(q^{\binom{m-t-1+k}k})\cr
=\,&\sum_{i=0}^t\alpha_i N_k(m-i)+O(q^{\binom{m-t-1+k}k}).
\end{split}
\]
\end{proof}

When $t=0$ and $1$ in \eqref{9A}, we obtain the asymptotic formulas 
in \cite{Car63} and \cite{Bod}. When $t=2$, equation~\eqref{9A} becomes 
\[
\begin{split}
I_k(m)=\,&N_k(m)-\frac{q(q^k-1)}{q-1}N_k(m-1)+\Bigl[\frac{q^2(q^k-1)^2}{(q-1)^2}-\frac{q^{k+1}(q^{\frac 12k(k+1)}-1)}{q-1}\Bigr]N_k(m-2)\cr
&+O(q^{\binom{m-3+k}k}).
\end{split}
\]


\section{An algorithm for calculating GCDs of polynomials in several variables}

In one variable, the Euclidean Algorithm can be used
to calculate the greatest common divisor of two polynomials. An analogous
algorithm for computing the greatest common divisor of polynomials in several variables is not so well known.
We next describe such an algorithm based on the arithmetic of the 
polynomial ring over a unique factorization domain and an induction on the number of variables.
(There is another algorithm for computing the greatest common divisor of polynomials in several variables using Gr\"obner bases, see \cite[Example 2.3.8]{Ada94}.)

Let $F$ be a field and let
\begin{gather*}
f(x,y)=f_0(y)+\cdots+f_m(y)x^m,\\
g(x,y)=g_0(y)+\cdots+g_n(y)x^n
\end{gather*}
be two polynomials in $F[x,y]$, where $y=(y_1,\dots,y_k)$ and 
$f_m,g_n\ne0$. To find $\text{gcd}_{F[x,y]}(f,g)$, we may assume $\text{gcd}_{F[y]}(f_0,\dots,f_m)=1$
and $\text{gcd}_{F[y]}(g_0,\dots,g_n)=1$, i.e., $f,g$ are primitive polynomials in $(F[y])[x]$. Then $\text{gcd}_{F[x,y]}(f,g)$ is the gcd of
$f$ and $g$ in $(F(y))[x]$ which belongs to $(F[y])[x]$ and is primitive. Therefore we have the following algorithm.

Let $h_0=f, h_1=g$ and rewrite
\[
h_i=h_i^{(0)}(y)+h_i^{(1)}(y)x+\cdots+h_i^{(d_i)}(y)x^{d_i},\qquad i=0,1,
\]
where $h_i^{(d_i)}\ne0$. 

Compute $h_2,h_3,\cdots\in (F[y])[x]$ inductively as follows. Switch $h_{i-1}$ and $h_i$ if necessary to make $d_i\le d_{i-1}$. Compute
\[
h_{i+1}=h_i^{(d_i)}h_{i-1}-h_{i-1}^{(d_{i-1})}x^{d_{i-1}-d_i}h_i
\]
and write the result as
\[
h_{i+1}=h_{i+1}^{(0)}(y)+h_{i+1}^{(1)}(y)x+\cdots+h_{i+1}^{(d_{i+1})}(y)x^{d_{i+1}}.
\]

Note that $\deg_x h_i$ is decreasing with respect to $i$.  
When it first occurs $h_{I+1}=0$, we have $h_I=\text{gcd}_{(F(y))[x]}(f,g)$. Hence
\[
\text{gcd}_{F[x,y]}(f,g)=\frac 1{\text{gcd}_{F[y]}(h_I^{(0)}(y),\dots,h_I^{(d_I)}(y))}h_I(x,y),
\]
where the denominator is the greatest common divisor in $k$ variables. The algorithm proceeds with induction on the number of variables.


\section{The number of  relatively prime pairs} 


\subsection{Formula for the number of  relatively prime pairs}

In this subsection we establish formulas for the number $P_k(m;n)$ of pairs  of normalized polynomials 
in $k$ variables of total degrees $m$ and $n$ over $\Bbb F_q$ which have greatest common divisor
$1$. The formula depends on $I_k(d)$, $1\le d\le \min\{m,n\}$.

Let $m,n\ge 0$. For $h\in \mathcal N_k(d)$ with $d\le\min\{m,n\}$, let 
\begin{gather*}
N_{=}(h)=|\{(f,g)\in \mathcal N_k(m)\times \mathcal N_k(n):\text{gcd}(f,g)=h\}|,\\
N_{\ge}(h)=|\{(f,g)\in \mathcal N_k(m)\times \mathcal N_k(n):h\mid \text{gcd}(f,g)\}|.
\end{gather*}
Then
\begin{equation}\label{Nge}
N_{\ge}(h)=N_k(m-d)N_k(n-d).
\end{equation}  
Since $N_{\ge}(h)=\sum_{h\mid u}N_{=}(u)$, by M\"obius inversion we have
\begin{equation}\label{Pkmn}
P_k(m;n)=N_{=}(1)=\sum_{h:\, 0\le\deg h\le\min\{m,n\}}\mu(h)N_{\ge}(h),
\end{equation}
where $\mu$ is the M\"obius function of the partially ordered set $(\mathcal N_k,\;\mid\;)$ and is given by
\[
\mu(h)=\begin{cases}
(-1)^s&\text{if $h$ is a product of $s$ distict irreducibles},\cr
0&\text{if $h$ is divisible by the square of an irreducible}.
\end{cases}
\]
By \eqref{Pkmn} and \eqref{Nge}, we have
\begin{equation}\label{2}
P_k(m;n)=\sum_{0\le d\le\min\{m,n\}}N_k(m-d)N_k(n-d)A_k(d),
\end{equation}
where
\[
A_k(d)=\sum_{h\in \mathcal N_k(d)}\mu(h).
\]
If $h\in\mathcal N_k(d)$ is such that $\mu(h)\ne 0$, then $h$ has to be a product of $a_1+\cdots+a_d$ distinct irreducibles, $a_i$ of which have degree $i$ ($1\le i\le d$), for some $a_1,\dots,a_d\in\Bbb N$ with $1a_1+2a_2+\cdots+da_d=d$; in such case, $\mu(h)=(-1)^{a_1+\cdots+a_d}$. Therefore 
\begin{equation}\label{Ad}
A_k(d)=\sum_{1a_1+2a_2+\cdots+da_d=d}(-1)^{a_1+\cdots+a_d}\binom{I_k(1)}{a_1}\cdots\binom{I_k(d)}{a_d}.
\end{equation}
In \eqref{2},
\[
N_k(m-d)N_k(n-d)=q^{\binom{m-d+k-1}k+\binom{n-d+k-1}k}\frac{(q^{\binom{m-d+k-1}{k-1}}-1)(q^{\binom{n-d+k-1}{k-1}}-1)}{(q-1)^2}.
\]
In \eqref{Ad}, $I_k(1),\dots,I_k(d)$ can be computed inductively by \eqref{2a}.

\medskip

\noindent{\bf Remark.}
When $k\ge 2$, no closed formula for $A_k(d)$ is known. When $k=1$, Carltz \cite{Car32} determined that
\begin{equation}\label{A1d}
A_1(d)=
\begin{cases}
1&\text{if}\ d=0,\cr
-q&\text{if}\ d=1,\cr
0&\text{if}\ d\ge 2.
\end{cases}
\end{equation}
Equations \eqref{2} and \eqref{A1d} provide yet another quick determination for $P_1(m;n)$ (cf. \cite{BenBen, Cor, Rei}).  

\medskip

We now consider the situation where $n$ is small and fixed and 
$m$ ($\ge n$) is arbitrary. Then $A_k(0),\dots,A_k(n)$, and hence formula \eqref{2}, can be made explicit. The first few terms of the sequence $A_k(d)$ are given below. 
\begin{gather*}
A_k(0)=1,\\
A_k(1)=-\frac{q(q^k-1)}{q-1},\\
A_k(2)=\frac q{(q-1)^2}\Bigl[q-q^k-q^{k+1}+q^{2k+1}-q^{\frac 12k(k+3)}(q-1)\Bigr],\\
\begin{split}
A_k(3)=\,&\frac{-1}{(q-1)^3}\Bigl[(q^k-1)(q^3-2q^{k+2}+q^{2k+3}+2q^{2+\frac 12k(k+3)}-2q^{3+\frac 12k(k+3)})\cr
&+q^{\frac 12(k+1)(k+2)}(q-1)^2(q^{\frac 16k(k+1)(k+2)}-1)\Bigr].
\end{split}
\end{gather*}


\subsection{Asymptotic results}

Here we prove some asymptotic results concerning the number 
${P_k(m;n)}$. When $k=1$, Theorem~\ref{T1.3} (with $s=2$) states that 
$\frac {P_1(m;n)}{N_1(m)N_1(n)}=1-\frac 1q$. When $k\ge 2$ the 
situation is totally different as shown in the next theorem. What causes 
this fundamental difference is the fact that almost all polynomials in one 
variable are 
reducible but almost all polynomials in more than one variable are 
irreducible. We will use the fact that when $k\ge 2$, $\frac{I_k(m)}{N_k(m)}\to 1$ 
as $m\to \infty$, which was established in \cite{Car63} and of course also 
follows from Theorem~\ref{T2.2}. 

\begin{Thm}\label{T4.1}
Let $k\ge 2$.  Then
\[
\lim_{m+n\to\infty}\frac{P_k(m;n)}{N_k(m)N_k(n)}=1.
\]
\end{Thm}

\begin{proof}
Without loss of generality we may assume $m\ge n$. Then
\[
\begin{split}
0\,&\le 1-\frac{P_k(m;n)}{N_k(m)N_k(n)}\cr
&=\frac{|\{(f,g)\in \mathcal N_k(m)\times \mathcal N_k(n):\text{gcd}(f,g)\ne 1\}|}{N_k(m)N_k(n)}\cr
&\le \frac{|\{(f,g)\in \mathcal N_k(m)\times \mathcal N_k(n):f\notin \mathcal I_k(m)\}|+|\{(f,f)\in \mathcal N_k(m)\times\mathcal  N_k(n):f\in \mathcal I_k(m)\}|}{N_k(m)N_k(n)}\cr
&\kern4cm (|\{(f,f)\in \mathcal N_k(m)\times\mathcal  N_k(n):f\in \mathcal I_k(m)\}|=0\ \text{if}\ m\ne n)\cr
&\le\frac{N_k(m)-I_k(m)+1}{N_k(m)}\to 0 \qquad\text{as}\ m\to\infty.
\end{split}
\]
\end{proof}

For values of $P_2(m;n)$ versus $N_2(m)N_2(n)$ with $q=2$ and $m,n\le 5$, see Table~\ref{Tb2} in Appendix B.

\medskip

Theorem~\ref{T4.1} can be restated as $P_k(m;n)=N_k(m)N_k(n)+
o(N_k(m)N_k(n))$ as $m+n\to \infty$. The following theorem gives an asymptotic formula for $P_k(m;n)$ up to $O(N_k(m-t-1)N_k(n-t-1))$ for any fixed $t\ge 0$.

\begin{Thm}\label{T4.2}
Let $k\ge 2$ and $t\ge 0$ be fixed integers. Then
\begin{equation}\label{eq-T4.2}
P_k(m;n)=\sum_{d=0}^t N_k(m-d)N_k(n-d)A_k(d)+O(N_k(m-t-1)N_k(n-t-1)),
\end{equation}
where $A_k(d)$ is defined in \eqref{Ad}. The constant in the $O$-term depends only on $q,k,t$.
\end{Thm} 

\begin{Lem}\label{L4.3}
Let $k\ge 2$ and $t\ge 0$. There exist constants $\epsilon=\epsilon(k,t)>0$ and $N=N(k,t)>0$ such that when $m+n\ge N$ and $1\le d\le \min\{m,n\}-t$,
\[
\binom{m-d+k}k+\binom{n-d+k}k+\binom{d+t+k}k\le \binom{m+k}k+\binom{n+k}k-\epsilon\cdot(m+n).
\]
\end{Lem}

\begin{proof}
Use induction on $k+t$. First assume $k=2$. Choose $N(2,t)>0$ such that $N(2,t)\ge 4(t^2+5t-1)$ and assume $m+n\ge N(2,t)$. We have

\[
\begin{split}
&2\Bigl[\binom{m-d+2}2+\binom{n-d+2}2+\binom{d+t+2}2\Bigr]\cr
=\,&(m-d+1)(m-d+2)+(n-d+1)(n-d+2)+(d+t+1)(d+t+2)\cr
=\,&(m+1)(m+2)-d(2m+3)+d^2+(n+1)(n+2)-d(2n+3)+d^2\cr
&+d^2+d(2t+3)+(t+1)(t+2)\cr
\le\,&(m+1)(m+2)+(n+1)(n+2)-d[2(m+n+3)-3d-(2t+3)-(t+1)(t+2)]\cr
=\,&(m+1)(m+2)+(n+1)(n+2)-d[\frac32(m+n)-3d+\frac 12(m+n)-(t^2+5t-1)]\cr
\le\,&(m+1)(m+2)+(n+1)(n+2)-d[\frac 12(m+n)-(t^2+5t-1)]\cr
\le\,&(m+1)(m+2)+(n+1)(n+2)-\frac 14(m+n)\qquad (\because m+n\ge 4(t^2+5t-1))\cr
=\,&2\Bigl[\binom{m+2}2+\binom{n+2}2-\frac 18(m+n)\Bigr].
\end{split}
\]

Now assume $k>2$ and $t=0$. Let $N(k,0)=N(k-1,0)$ and assume $m+n\ge N(k,0)$. We have
\begin{equation}\label{13a}
\begin{split}
&\binom{m-d+k}k+\binom{n-d+k}k+\binom{d+k}k\cr
=\,&\binom{m-d+k-1}{k-1}+\binom{n-d+k-1}{k-1}+\binom{d+k-1}{k-1}\cr
&+\binom{m-d+k-1}k+\binom{n-d+k-1}k+\binom{d+k-1}k\cr
\le\,&\binom{m+k-1}{k-1}+\binom{n+k-1}{k-1}-\epsilon(k-1,0)(m+n)\cr
&+\binom{m-d+k-1}k+\binom{n-d+k-1}k+\binom{d+k-1}k\cr
&\kern4cm\text{(since $m+n\ge N(k-1,0)$, the induction hypothesis applies)}.
\end{split}
\end{equation}
Note that 
\begin{equation}\label{14a}
\binom{n-d+k-1}k+\binom{d+k-1}k\le\binom{n+k-1}k.
\end{equation}
The quickest way to see \eqref{14a} is to observe that $\binom{n-d+k-1}k$,
 respectively, $\binom{d+k-1}k$, $\binom{n+k-1}k$, is the number of ways to ``choose $k$ from $n-d$,  respectively, $d$, $n$, with repetition''.
Thus \eqref{13a} continues as
\[
\begin{split}
&\binom{m-d+k}k+\binom{n-d+k}k+\binom{d+k}k\cr
\le\,&\binom{m+k-1}{k-1}+\binom{m+k-1}k+\binom{n+k-1}{k-1}+\binom{n+k-1}k-\epsilon(k-1,0)(m+n)\cr
=\,&\binom{m+k}k+\binom{n+k}k-\epsilon(k-1,0)(m+n).
\end{split}
\]

Now assume $k>2$ and $t>0$. Let $N(k,t)=\max\{N(k-1,t), N(k,t-1)+2\}$ and assume $m+n\ge N(k,t)$.  We have

\[
\begin{split}
&\binom{m-d+k}k+\binom{n-d+k}k+\binom{d+t+k}k\cr
=\,&\binom{m-d+k-1}{k-1}+\binom{n-d+k-1}{k-1}+\binom{d+t+k-1}{k-1}\cr
&+\binom{m-1-d+k}k+\binom{n-1-d+k}k+\binom{d+t-1+k}k\cr
\le\,&\binom{m+k-1}{k-1}+\binom{n+k-1}{k-1}-\epsilon(k-1,t)(m+n)\cr
&+\binom{m-1+k}k+\binom{n-1+k}k-\epsilon(k,t-1)(m-1+n-1)\cr
&\text{(since $m+n\ge N(k-1,t)$ and $m-1+n-1\ge N(k,t-1)$, the induction hypothesis applies)}\cr
\le\,&\binom{m+k}k+\binom{n+k}k-\epsilon(k-1,t)(m+n).
\end{split}
\]
\end{proof}

\begin{proof}[Proof of Theorem~\ref{T4.2}]
When $\min\{m,n\}\le t$, the $O$-term in \eqref{eq-T4.2} is $0$. So we assume $m,n>t$. We have by \eqref{2}
\[
P_k(m;n)=\sum_{d=0}^t N_k(m-d)N_k(n-d)A_k(d)+\sum_{d=t+1}^{\min\{m,n\}}N_k(m-d)N_k(n-d)A_k(d),
\]
where
\[
\begin{split}
\Bigl|\sum_{d=t+1}^{\min\{m,n\}}N_k(m-d)N_k(n-d)A_k(d)\Bigr|
\le\,&\sum_{d=t+1}^{\min\{m,n\}}N_k(m-d)N_k(n-d)|A_k(d)|\cr
\le\,&\sum_{d=t+1}^{\min\{m,n\}}N_k(m-d)N_k(n-d)N_k(d).
\end{split}
\]
So it suffices to show
\begin{equation}\label{16a}
\sum_{d=t+1}^{\min\{m,n\}}N_k(m-d)N_k(n-d)N_k(d)=O(N_k(m-t-1)N_k(n-t-1)).
\end{equation}
We have 
\begin{equation}\label{17a}
\begin{split}
&\sum_{d=t+1}^{\min\{m,n\}}N_k(m-d)N_k(n-d)N_k(d)\cr
=\,&N_k(m-t-1)N_k(n-t-1)N_k(t+1)+\sum_{d=t+2}^{\min\{m,n\}}N_k(m-d)N_k(n-d)N_k(d)\cr
\le\,&O(N_k(m-t-1)N_k(n-t-1))+\sum_{d=t+2}^{\min\{m,n\}} q^{\binom{m-d+k}k+\binom{n-d+k}k+\binom{d+k}k}.
\end{split}
\end{equation}
By Lemma~\ref{L4.3}, for $t+2\le d\le\min\{m,n\}$,
\[
\begin{split}
&\binom{m-d+k}k+\binom{n-d+k}k+\binom{d+k}k\cr
=\,&\binom{m-t-1-(d-t-1)+k}k+\binom{n-t-1-(d-t-1)+k}k+\binom{(d-t-1)+t+1+k}k\cr
\le\,&\binom{m-t-1+k}k+\binom{n-t-1+k}k-\epsilon(k,t+1)(m+n-2t-2)\cr
&\kern4cm \text{(when $m+n-2t-2\ge N(k,t+1)$)}.
\end{split}
\]
So
\[
\begin{split}
\sum_{d=t+2}^{\min\{m,n\}} q^{\binom{m-d+k}k+\binom{n-d+k}k+\binom{d+k}k}
\le\,&(m+n)q^{\binom{m-t-1+k}k+\binom{n-t-1+k}k-\epsilon(k,t+1)(m+n-2t-2)}\cr
=\,&o(1)\cdot q^{\binom{m-t-1+k}k+\binom{n-t-1+k}k}\cr
=\,&o(1)\cdot O(N_k(m-t-1)N_k(n-t-1))\cr
=\,&o(N_k(m-t-1)N_k(n-t-1)).
\end{split}
\]
Combining this with \eqref{17a} we arrive at \eqref{16a}.
\end{proof}


\section{Corresponding results for the vector degree} 

In Sections 2 and 4 we have considered the total degree versions of 
Problems~\ref{Pb1.1} and \ref{Pb1.2}. In this section we turn to the vector degree versions of the problems.
We will see that results similar to those in Sections 2 and 4 also hold for the vector degree. However, the proofs are not always simple parallels of those  in the total degree case. In fact, asymptotic results in the vector degree case are considerably more difficult to prove than in the total degree case.


\subsection{Recursive formula for $I_k(\frak m)$}

Recall that for $\frak m=(m_1,\dots,m_k)\in\Bbb N^k$, $N_k(\frak m)=|\mathcal N_k(\frak m)|$ where $\mathcal N_k(\frak m)=\{f\in \mathcal N_k:\text{Deg}\,f=\frak m\}$ and 
$I_k(\frak m)=|\mathcal I_k(\frak m)|$ where $\mathcal I_k(\frak m)=\{f\in \mathcal N_k(\frak m):\text{$f$ is irreducible}\}$. 
For $\frak m=(m_1,\dots,m_k), \frak n=(n_1,\dots,n_k)\in\Bbb N^k$, we write $\frak n\le \frak m$ to mean that $n_i\le m_i$ for all $1\le i\le k$;
$\frak n<\frak m$ means that $\frak n\le \frak m$ and $\frak n\ne \frak m$. The zero tuple $(0,\dots,0)\in\Bbb N^k$ is denoted by $\frak 0$.  
Define   
\[
\mathcal N_\le(\frak m)=\{0\ne f\in \mathcal N_k:\ \text{Deg}\,f\le\frak m\}.
\]
Then
\[
\sum_{\frak n\le \frak m}N_k(\frak n)=|\mathcal N_\le(\frak m)|=\frac{q^{(m_1+1)\cdots(m_k+1)}-1}{q-1}.
\]
The M\"obius function of $(\Bbb N^k,\le)$ is
\[
\mu(\frak m,\frak n)=
\begin{cases}
(-1)^{m_1+\cdots+m_k-n_1-\cdots-n_k}&\text{if}\ \frak m-\frak n\in\{0,1\}^k,\cr
0&\text{otherwise}.
\end{cases}
\]
So by M\"obius inversion,

\[
\begin{split}
N_k(\frak m)\,&=\sum_{\frak n\le \frak m}\mu(\frak m,\frak n)\frac{q^{(n_1+1)\cdots(n_k+1)}-1}{q-1}\cr
&=\sum_{(\delta_1,\dots,\delta_k)\in\{0,1\}^k}(-1)^{\delta_1+\dots+\delta_k}\frac{q^{(m_1-\delta_1+1)\cdots(m_k-\delta_k+1)}-1}{q-1}\cr
&=\frac 1{q-1}\sum_{(\delta_1,\dots,\delta_k)\in\{0,1\}^k}(-1)^{k+\delta_1+\dots+\delta_k}(q^{(m_1+\delta_1)\cdots(m_k+\delta_k)}-1)\cr
&=\frac 1{q-1}\sum_{(\delta_1,\dots,\delta_k)\in\{0,1\}^k}(-1)^{k+\delta_1+\dots+\delta_k}q^{(m_1+\delta_1)\cdots(m_k+\delta_k)}.
\end{split}
\]
This formula is (3.1) in Cohen \cite{Coh}. Our proof is different from that of \cite{Coh}.

Unique factorization in $\Bbb F_q[x_1,\dots,x_k]$ gives
\begin{equation}\label{11}
\sum_{\substack{(a_{\frak i})_{\frak 0<\frak i\le\frak m}\cr \sum_{\frak i}a_{\frak i}\frak i=\frak m}}\prod_{\frak i}\binom{I_k(\frak i)+
a_{\frak i}-1}{a_{\frak i}}=N_k(\frak m).
\end{equation}
Hence $I_k(\frak m)$ can be obtained  inductively by $I_k(\frak 0)=0$ and
\begin{equation}\label{11b}
I_k(\frak m)=N_k(\frak m)-\sum_{\substack{(a_{\frak i})_{\frak 0<\frak i<\frak m}\cr \sum_{\frak i}a_{\frak i}\frak i=\frak m}}\prod_{\frak i}\binom{I_k(\frak i)+
a_{\frak i}-1}{a_{\frak i}},\qquad \frak m>\frak 0.
\end{equation}
A table of $I_2(m_1,m_2)$ with $q=2$ and $0\le m_1\le m_2\le 5$ is given in Appendix B.

Cohen \cite{Coh} computed $I_k(\frak m,n)$ explicitly for 
$\frak m=(1),(2),(3),(1,1),(1,2)$ and arbitrary $n\in \Bbb N$. In general if 
$\frak m\in\Bbb N^{k-1}$ is small, an explicit formula for $I_k(\frak m,n)$ 
can be obtained using \eqref{11} and induction. To illustrate the method, we 
include the computation for $I_2(1,n)$ and $I_2(2,n)$. 

When $\frak m=(1,n)$, \eqref{11} can be written as
\[
N_2(1,n)=\sum_{0\le t\le n}I_2(1,t)N_1(n-t)=\sum_{0\le t\le n}I_2(1,t)q^{n-t}.
\]
So
\[
\sum_{0\le t\le n}q^{-t}I_2(1,t)=q^{-n}N_2(1,n).
\]
It follows that for $n>0$
\[
q^{-n}I_2(1,n)=q^{-n}N_2(1,n)-q^{-(n-1)}N_2(1,n-1).
\]
Thus
\begin{equation}\label{14b}
I_2(1,n)=N_2(1,n)-qN_2(1,n-1)=(q^2-1)q^{2n-1}.
\end{equation}

When $\frak m=(2,n)$, \eqref{11} can be written as
\[
\begin{split}
&N_2(2,n)\cr
=\,&\sum_{0\le t\le\frac n2}\binom{I_2(1,t)+1}2N_1(n-2t)\cr
&+\sum_{\substack{0\le s<t\cr s+t\le n}}I_2(1,s)I_2(1,t) N_1(n-s-t)
+\sum_{0\le t\le n}I_2(2,t)N_1(n-t)\cr
=\,&\sum_{0\le t\le\frac n2}\binom{I_2(1,t)+1}2 q^{n-2t}+\sum_{\substack{0\le s<t\cr s+t\le n}}I_2(1,s)I_2(1,t) q^{n-s-t}
+\sum_{0\le t\le n}I_2(2,t)q^{n-t},
\end{split}
\]
namely
\[
\sum_{0\le t\le n}q^{-t}I_2(2,t)
=q^{-n}N_2(2,n) -\sum_{0\le t\le\frac n2}\binom{I_2(1,t)+1}2 q^{-2t}-\sum_{\substack{0\le s<t\cr s+t\le n}}I_2(1,s)I_2(1,t) q^{-s-t}.
\]
It follows that for $n>0$
\[
\begin{split}
q^{-n}I_2(2,n)=\,&q^{-n}N_2(2,n)-q^{-(n-1)}N_2(2,n-1)-\delta(\frac n2)\binom{I_2(1,\frac n2)+1}2 q^{-n}\cr
&-\sum_{\substack{0\le s<t\cr s+t= n}}I_2(1,s)I_2(1,t) q^{-n},
\end{split}
\]
where
\[
\delta(x)=
\begin{cases}
1&\text{if}\ x\in\Bbb Z,\cr
0&\text{otherwise}.
\end{cases}
\] 
So
\begin{equation}\label{15b}
\begin{split}
&I_2(2,n)\cr
=\,&N_2(2,n)-qN_2(2,n-1)-\delta(\frac n2)\binom{I_2(1,\frac n2)+1}2+\frac 12\delta(\frac n2)I_2(1,\frac n2)^2-\frac 12\sum_{s+t=n}I_2(1,s)I_2(1,t)\cr
=\,&N_2(2,n)-qN_2(2,n-1)-\delta(\frac n2)\frac12 I_2(1,\frac n2)-\frac 12 \sum_{\substack{s+t=n\cr s,t>0}}(q^2-1)^2q^{2(s+t)-2}-I_2(1,0)I_2(1,n)\cr
=\,&(q^2-1)\Bigl[q^{3n-2}(1+q+q^2)-\frac 12 q^{2n-2}\bigl(n(q^2-1)+(q+1)^2\bigr)-\frac 12\delta(\frac n2)q^{n-1}\Bigr].
\end{split}
\end{equation}

\noindent
{\bf Remark.}
\begin{itemize}
\item[(i)] Recall that $I_2(0,n)=I_1(n)=\frac 1n\sum_{d\mid n}\mu(\frac nd)q^d$. In comparison, the formulas for $I_2(1,n)$, $I_2(2,n),\dots$ are more explicit 
and do not involve the M\"obius function. 

\item[(ii)]
For fixed (and small) $m>0$ and arbitrary $\frak n\in\Bbb N^{k-1}$ where $k\ge 3$, it is not clear what kind of formula one might expect for $I_k(m,\frak n)$. 
\end{itemize}


\subsection{Asymptotic formula for $I_k(\frak m)$}

When $k\ge 2$, the notion $(m_1,\dots,m_k)\to\infty$ is rather ambiguous; its 
precise meaning depends on further assumptions on $m_1,\dots,m_k$.
For example, in some theorems of \cite{Coh}, $m_1,\dots,m_{k-1}$ are fixed and $m_k$ is allowed to approach $\infty$. In our discussion of asymptotic formula for $I_k(m_1,\dots,m_k)$, we  will not fix any of $m_1,\dots,m_k$. The notation $O(A(m_1,\dots,m_k))$ always denotes a function $B(m_1,\dots,m_k)$ such that 
\[
|B(m_1,\dots,m_k)|\le CA(m_1,\dots,m_k)
\]
for all $m_1,\dots, m_k$ in the described range, where $C>0$ is a constant independent of $m_1,\dots,m_k$.

We will need the following key lemma which plays the same role as
 Lemma~\ref{L2.1} played  in the total degree case.

\begin{Lem}\label{L0.3}
Let $k\ge 2$ and
let $\frak t=(t_1,0,\dots,0)\in\Bbb N^k$ be fixed. Let $\frak m=(m_1,\dots,m_k)\in(\Bbb Z^+)^k$. Assume $m_1>2t_1$ when $k=2$ and $m_1\ge 2t_1$ when $k\ge 3$. Then
\begin{equation}\label{16b}
\sum_{\substack{(a_{\frak i})_{\frak 0<\frak i\le \frak m-\frak t}\cr
\sum_{\frak i}a_{\frak i}\frak i=\frak m}}\prod_{\frak 0<\frak i\le \frak m-\frak t}
\binom{I_k(\frak i)+a_{\frak i}-1}{a_{\frak i}}=O(q^{(m_1-t_1+1)(m_2+1)\cdots(m_k+1)}).
\end{equation}
\end{Lem}

The sum in \eqref{16b} is the number of $f\in\mathcal N_k(\frak m)$ whose irreducible factors all have $\deg_{x_1}\le m_1-t_1$.
Lemma~\ref{L0.3} is not difficult to see from an intuitive point of view, but its proof is rather tedious. The proof of Lemma~\ref{L0.3} is given in Appendix A.

\medskip

\noindent {\bf Note.} 
If $k=2$ and $t_1<m_1\le 2t_1$, \eqref{16b} does not hold. In fact, since $m_1-t_1\le\frac 12 m_1$, the sum in \eqref{16b} is 
\[
\ge \sum_{\frac{m_2}4\le d<\frac {m_2}2}I_2(m_1-t_1,d)I_2(m_1-t_1,m_2-d).
\]
If $m_1-t_1\ge 3$, from Theorem~\ref{T5.3} below we can derive that for $d>0$,
\[
I_2(m_1-t_1,d)=q^{(m_1-t_1+1)(d+1)}(1-q^{-m_1-t_1})(1-q^{-m_1-t_1-1})+O(q^{(m_1-t_1)(d+1)}).
\]
Hence for $\frac{m_2}4\le d<\frac{m_2}2$, 
\[
\begin{split}
I_2(m_1-t_1,d)I_2(m_1-t_1,m_2-d)=\,&q^{(m_1-t_1+1)(m_2+2)}(1-q^{-m_1-t_1})(1-q^{-m_1-t_1-1})\cr
&+O(q^{(m_1-t_1+1)(m_2+2)-\frac{m_2}4}).
\end{split}
\]
So
\[
\begin{split}
&\sum_{\frac{m_2}4\le d<\frac {m_2}2}I_2(m_1-t_1,d)I_2(m_1-t_1,m_2-d)\cr
=\,&\left\lfloor\frac{m_2}4\right\rfloor q^{(m_1-t_1+1)(m_2+2)}(1-q^{-m_1-t_1})(1-q^{-m_1-t_1-1})+O(m_2q^{(m_1-t_1+1)(m_2+2)-\frac{m_2}4})\cr
\ne\,&O(q^{(m_1-t_1+1)(m_2+1)}).
\end{split}
\]
If $m_1-t_1=1$ or $2$, using \eqref{14b} and \eqref{15b}, similar arguments show that \eqref{16b} also fails.

\begin{Lem}\label{LB.3}
Let $k\ge 2$, $m_1,\dots,m_k\in \Bbb Z^+$ and $a,b\in\Bbb N$. Assume $m_1=\max_{1\le i\le k}m_i$ and let $1\le s<k$. Then 
\[
\begin{split}
&(a+1)\prod_{i=1}^s(m_i+1)+(b+1)\prod_{i=s+1}^k(m_i+1)\cr
\le\,&
\begin{cases}
\displaystyle (a+b+1) m_1\prod_{i=2}^k(m_i+1)-a-b+1&\text{if}\ m_1\ge 2, \vspace{2mm}\cr 
\displaystyle (a+b+1)m_1\prod_{i=2}^k(m_i+1)-\min\{a,b\}+2&\text{if $m_1=1$ and $k\ge 3$}.
\end{cases}
\end{split}
\] 
\end{Lem}

\begin{proof}
First assume $m_1\ge 2$.
Let $m_{i_0}=\min_{1\le i\le k}m_i$. We claim that 
\begin{equation}\label{23a}
(a+1)\prod_{i=1}^s(m_i+1)+(b+1)\prod_{i=s+1}^k(m_i+1)\le
(a+b+1)\prod_{i\ne i_0}(m_i+1)+m_{i_0}+1.
\end{equation}
Without loss of generality, assume $\prod_{i=1}^s(m_i+1)\ge\prod_{i=s+1}^k(m_i+1)$. Then

\[
\begin{split}
&(a+1)\prod_{i=1}^s(m_i+1)+(b+1)\prod_{i=s+1}^k(m_i+1)\cr
\le\,&(a+b+1)\prod_{i=1}^s(m_i+1)+\prod_{i=s+1}^k(m_i+1)\cr
=\,&(a+b+1)\prod_{i=1}^s(m_i+1)+\Bigl[\prod_{i=s+1}^{k-1}(m_i+1)\Bigr](m_k+1)\cr
\le\,&\Bigl[(a+b+1)\prod_{i=1}^s(m_i+1)\Bigr]\prod_{i=s+1}^{k-1}(m_i+1)+m_k+1\qquad(\because (a+b+1)\prod_{i=1}^s(m_i+1)\ge m_k+1)\cr
=\,&(a+b+1)\prod_{i=1}^{k-1}(m_i+1)+m_k+1.
\end{split}
\]
The last expression does not decrease when interchanging $m_k$ and $m_{i_0}$. So \eqref{23a} holds.

For simplicity, assume $i_0=k$ in \eqref{23a}. We have
\[
\begin{split}
&(a+b+1)\prod_{i=1}^{k-1}(m_i+1)+m_k+1-(a+b+1)m_1\prod_{i=2}^k(m_i+1)\cr
=\,&[m_1+1-m_1(m_k+1)](a+b+1)\prod_{i=2}^{k-1}(m_i+1)+m_k+1\cr
=\,&-(m_1m_k-1)(a+b+1)\prod_{i=2}^{k-1}(m_i+1)+m_k+1\cr
\le\,&-[(m_1m_k-1)+(a+b+1)-1]+m_k+1\qquad(\because m_1m_k-1\ge 1)\cr
=\,&-a-b-m_1m_k+m_k+2\cr
\le\,&-a-b+1.
\end{split}
\]

Now assume $m_1=1$ and $k\ge 3$. Thus $m_1=\cdots=m_k=1$. Then
\[
\begin{split}
&(a+1)\prod_{i=1}^s(m_i+1)+(b+1)\prod_{i=s+1}^k(m_i+1)-(a+b+1)m_1\prod_{i=2}^k(m_i+1)\cr
=\,&2^s(a+1)+2^{k-s}(b+1)-2^{k-1}(a+b+1)\cr
\le\,&2^{k-1}(\max\{a,b\}+1)+2(\min\{a,b\}+1)-2^{k-1}(a+b+1)\cr
=\,&-2^{k-1}\min\{a,b\}+2(\min\{a,b\}+1)\cr
\le\,&-\min\{a,b\}+2.
\end{split}
\]
\end{proof}

For $\frak m=(m_1,\dots,m_k)\in\Bbb N^k$, $\frak m^{(j)}=(m^{(j)}_1,\dots,m^{(j)}_k)\in \Bbb N^k$, $1\le j\le s$, we write 
$\frak m=\frak m^{(1)}\oplus\cdots\oplus\frak m^{(s)}$ if
$\frak m=\frak m^{(1)}+\cdots+\frak m^{(s)}$ and  $m_i^{(j)}=0$ or $m_i$ for each $1\le i\le k$ and $1\le j\le s$.

\begin{Lem}\label{L5.2}
Let $k\ge 2$, $n\in\Bbb N$ and $\frak m=(m_1,\dots,m_k)\in(\Bbb Z^+)^k$ with $m_1=\max_{1\le i\le k}m_i$. Assume that $m_1\ge 2$ when $k=2$. Let $s\ge 2$ and let 
$\frak m^{(j)}=(m^{(j)}_1,\dots,m^{(j)}_k)\in\Bbb N^k$, $1\le j\le s$, such that 
$\frak m=\frak m^{(1)}\oplus\cdots\oplus\frak m^{(s)}$.
Then
\begin{equation}\label{18c}
\begin{split}
&|\{f\in \mathcal N_{k+1}(\frak m,n):f=f_1\cdots f_sg,\ f_j\in \mathcal N_{k+1}(\frak m^{(j)},*),\ g\in \mathcal N_{k+1}(\frak 0,*)\}|\cr
=\,&O(q^{(n+1)m_1\prod_{i=2}^k(m_i+1)}).
\end{split}
\end{equation}
\end{Lem}

\begin{proof}
 We have 
\begin{equation}\label{19c}
\begin{split}
\text{LHS of \eqref{18c}}
\le\,&\sum_{a_1+\cdots+a_{s+1}=n}N_{k+1}(\frak m^{(1)},a_1)\cdots N_{k+1}(\frak m^{(s)},a_s)N_{k+1}(\frak 0,a_{s+1})\cr
\le\,&\sum_{a_1+\cdots+a_{s+1}=n}q^{(a_1+1)\prod_{i=1}^k(m_i^{(1)}+1)+\cdots+(a_s+1)\prod_{i=1}^k(m_i^{(s)}+1)+(a_{s+1}+1)}.
\end{split}
\end{equation}
Note that
\[
\begin{split}
&(a_1+1)\prod_{i=1}^k(m_i^{(1)}+1)+\cdots+(a_s+1)\prod_{i=1}^k(m_i^{(s)}+1)+(a_{s+1}+1)\cr
\le\,&(a_1+1)\prod_{i=1}^k(m_i^{(1)}+1)+(a_2+\cdots+a_{s+1}+1)\prod_{i=1}^k(m_i^{(2)}+\cdots+m_i^{(s)}+1)+s-1\cr
&\kern8cm\text{(by Lemma~\ref{LB.1})}\cr
\le\,&(a_1+\cdots+a_{s+1}+1)m_1\prod_{i=2}^k(m_i+1)-\min\{a_1,a_2+\cdots+a_{s+1}\}+2+s-1\cr
&\kern8cm\text{(by Lemma~\ref{LB.3})}\cr
\le\,&(n+1)m_1\prod_{i=2}^k(m_i+1)-\min\{a_1,n-a_1\}+k+1.
\end{split}
\]
In the above it is clear that $a_1$ can be replaced with any $a_l$, $1\le l\le s$. Thus
\[
\begin{split}
&(a_1+1)\prod_{i=1}^k(m_i^{(1)}+1)+\cdots+(a_s+1)\prod_{i=1}^k(m_i^{(s)}+1)+(a_{s+1}+1)\cr
\le\,&k+1+(n+1)m_1\prod_{i=2}^k(m_i+1)-\frac 1s\sum_{l=1}^s\min\{a_l,n-a_l\}\cr
\le\,&k+1+(n+1)m_1\prod_{i=2}^k(m_i+1)-\frac 1k\sum_{l=1}^s\min\{a_l,n-a_l\}.
\end{split}
\]
Returning to \eqref{19c}, we have
\[
\begin{split}
\text{LHS of \eqref{18c}}
\le\,&q^{k+1}q^{(n+1)m_1\prod_{i=2}^k(m_i+1)}\sum_{a_1+\cdots+a_{s+1}=n}q^{-\frac 1k\sum_{l=1}^s\min\{a_l,n-a_l\}}\cr
\le\,&q^{k+1}q^{(n+1)m_1\prod_{i=2}^k(m_i+1)}\sum_{a_1,\dots,a_s\le n}q^{-\frac 1k\sum_{l=1}^s\min\{a_l,n-a_l\}}\cr
=\,&O(q^{(n+1)m_1\prod_{i=2}^k(m_i+1)}).
\end{split}
\]
\end{proof}

\begin{Thm}\label{T5.3}
Let $k\ge 2$ and $(m_1,\dots,m_k)\in(\Bbb Z^+)^k$ with $m_1=\max_{1\le i\le k-1}m_i$. Assume that $m_1\ge 3$ when $k=2$ and that $m_1\ge 2$ when $k=3$. Then
\[
I_k(m_1,\dots,m_k)=N_k(m_1,\dots,m_k)-q N_k(m_1,\dots,m_{k-1},m_k-1)+O(q^{m_1(m_2+1)\cdots(m_k+1)}).
\]
\end{Thm}

\noindent{\bf Note.}
Theorem~\ref{T5.3} indicates that most of polynomials in $\mathcal N_k(m_1,\dots,m_k)$ that fail to be irreducible are of the form $(x_k+\alpha)f$ for some $\alpha\in\Bbb F_q$ and $f\in\mathcal N_k(m_1,\dots,m_{k-1},m_k-1)$. The asymptotic formula in Theorem~\ref{T5.3} is interesting only when $m_k\ge m_1$ since otherwise the $O$-term is bigger that the term $q N_k(m_1,\dots,m_{k-1},m_k-1)$.

\begin{proof}[Proof of Theorem~\ref{T5.3}]
Let $\frak m=(m_1,\dots,m_k)$ and $\frak m'=(m_1,\dots,m_{k-1})$. Write
\begin{equation}\label{N1-3}
\mathcal N_k(\frak m)=\mathcal N^{(1)}\overset{\cdot}\cup\mathcal N^{(2)}\overset{\cdot}\cup\mathcal N^{(3)}
\end{equation}
where
\begin{gather*}
\mathcal N^{(1)}=\{f\in \mathcal N_k(\frak m):f=f_1g,\ f_1\in \mathcal I_k(\frak m',a)\ \text{for some $0\le a\le m_k$}\},\\
\mathcal N^{(2)}=\{f\in \mathcal N_k(\frak m):f=f_1\cdots f_sg,\ s\ge 2,\ f_j\in \mathcal I_k(\frak m^{(j)},*),\ \frak m'=\frak m^{(1)}\oplus+\cdots
\oplus\frak m^{(s)},\ \frak m^{(j)}\ne \frak 0\},\\
\mathcal N^{(3)}=\{f\in \mathcal N_k(\frak m):\text{$f$ has a factor $f_1$ such that $0<{\deg}_{x_i}f_1<m_i$ for some $1\le i\le k-1$}\}.
\end{gather*}
First of all, we have
\begin{equation}\label{A}
|\mathcal N^{(1)}|=\sum_{a=0}^{m_k}I_k(\frak m',a)q^{m_k-a}.
\end{equation}
When $k=2$, $\mathcal N^{(2)}=\emptyset$; when $k\ge 3$, by Lemma~\ref{L5.2}, we have
\begin{equation}\label{B}
\begin{split}
|\mathcal N^{(2)}|\,&\le O(q^{m_1(m_2+1)\cdots(m_k+1)})\cdot(\text{the number of partitions of $\{1,\dots,k\}$})\cr 
&=O(q^{m_1(m_2+1)\cdots(m_k+1)}).
\end{split}
\end{equation}
We also claim that
\begin{equation}\label{C}
|\mathcal N^{(3)}|\le O(q^{m_1(m_2+1)\cdots(m_k+1)}).
\end{equation}
When $k\ge 3$, we have
\[
\begin{split}
|\mathcal N^{(3)}|\,&\le\sum_{1\le i\le k-1\, :\, m_i\ge 2}|\{f\in \mathcal N_k(\frak m):\text{all irreducible factors of $f$ have ${\deg}_{x_i}\le m_i-1$}\}|\cr
&=\sum_{1\le i\le k-1 \,:\, m_i\ge 2}O(q^{\frac {m_i}{m_i+1}(m_1+1)\cdots(m_k+1)})\qquad\text{(by Lemma~\ref{L0.3})}\cr
&=O(q^{m_1(m_2+1)\cdots(m_k+1)}).
\end{split}
\]
When $k=2$, the sum in the above has only one term with $i=1$. Since $m_1\ge 3$ by assumption, Lemma~\ref{L0.3} still applies. 
Combining \eqref{N1-3} -- \eqref{C}, we have
\[
N_k(\frak m)=\sum_{a=0}^{m_k}I_k(\frak m',a)q^{m_k-a}+O(q^{m_1(m_2+1)\cdots(m_k+1)}).
\]
Thus
\[
\sum_{a=0}^{m_k}I_k(\frak m',a)q^{-a}=N_k(\frak m)q^{-m_k}+q^{-m_k}O(q^{m_1(m_2+1)\cdots(m_k+1)}).
\]
It follows that
\[
I_k(\frak m',m_k)q^{-m_k}=N_k(\frak m)q^{-m_k}-N_k(\frak m',m_k-1)q^{-(m_k-1)}+q^{-m_k}O(q^{m_1(m_2+1)\cdots(m_k+1)}),
\]
i.e.,
\[
I_k(\frak m)=N_k(\frak m)-qN_k(\frak m',m_k-1)+O(q^{m_1(m_2+1)\cdots(m_k+1)}).
\]
\end{proof} 

Theorem~\ref{T5.3} improves the main result (Theorem 1) of \cite{Coh} by removing the factor $m_im_k$ in the $O$-term. It is also an improvement of Theorem 2 of \cite{Coh} since $m_1,\dots,m_{k-1}$ in Theorem~\ref{T5.3} are not fixed.

\medskip

\noindent{\bf Note.} Although we will not pursue further improvement of Theorem~\ref{T5.3} in the present paper, we mention that it is possible to improve the $O$-term in Theorem~\ref{T5.3} to $O(q^{(m_1-t_1+1)(m_2+1)\cdots(m_k+1)})$ for a fixed $t_1>0$.


\subsection{Number of relatively prime pairs}

Let $\frak m,\frak n\in\Bbb N^k$. For $h\in \mathcal N_k(\frak d)$ with $\frak d\le\frak m$ and $\frak d\le\frak n$, let 
\begin{gather*}
N_{=}(h)=|\{(f,g)\in \mathcal N_k(\frak m)\times\mathcal  N_k(\frak n):\text{gcd}(f,g)=h\}|,\\
N_{\ge}(h)=|\{(f,g)\in \mathcal N_k(\frak m)\times\mathcal  N_k(\frak n):h\mid \text{gcd}(f,g)\}|=\sum_{h\mid u}N_{=}(u).
\end{gather*}
Then
\[
N_{\ge}(h)=N_k(\frak m-\frak d)N_k(\frak n-\frak d).
\]  
By  M\"obius inversion we have
\begin{equation}\label{19a}
P_k(\frak m;\frak n)=N_{=}(1)=\sum_{h\,:\, \frak 0\le\text{Deg}\, h\le\frak m,\frak n}\mu(h)N_{\ge}(h)
=\sum_{\frak 0\le \frak d\le\frak m,\frak n}N_k(\frak m-\frak d)N_k(\frak n-\frak d)A_k(\frak d),
\end{equation}
where
\begin{equation}\label{A_k(d)}
A_k(\frak d)=\sum_{h\in\mathcal N_k(\frak d)}\mu(\frak d)=\sum_{\substack{(a_{\frak i})_{\frak 0< \frak i\le\frak d}\cr
\sum_{\frak 0< \frak i\le\frak d}a_{\frak i}\frak i=\frak d}}(-1)^{\sum_{\frak 0<\frak i\le \frak d}a_{\frak i}}\prod_{\frak 0< \frak i\le\frak d}
\binom{I_k(\frak i)}{a_{\frak i}},
\end{equation}
In \eqref{19a}, $I_k(\frak i)$ ($\frak 0< \frak i \le \frak d)$ can be computed inductively by \eqref{11b}.
The function $A_k(\frak d)$ in \eqref{A_k(d)} can be made explicit for small $\frak d$.

We include in Appendix B a table of $P_2(\frak m;\frak n)$ with $q=2$ and $\frak m,\frak n\le (4,4)$.


\subsection{Asymptotic formula for $P_k(\frak m;\frak n)$}

Let $k\ge 2$ and  $(m_1,\dots,m_k)\in\Bbb N^k$. It follows from Theorem~\ref{T5.3} that 
\[
\frac{I_k(m_1\dots,m_k)}{N_k(m_1,\dots,m_k)}\to 1\qquad \text{as $m_1\to\infty$ and $m_k\to\infty$}.
\]
By symmetry, the above statement holds as any two of the components of $(m_1,\dots,m_k)$ approach to $\infty$.

\begin{Thm}\label{T}
Let $k\ge 2$ and $\frak m=(m_1,\dots,m_k),\ \frak n=(n_1,\dots,n_k)\in\Bbb N^k$. Then
\[
\lim_{m_{k-1},m_k\to\infty}\frac{P_k(\frak m;\frak n)}{N_k(\frak m)N_k(\frak n)}=1.
\]
\end{Thm}

\begin{proof}
Let
\[
(\frak m';\frak n')=\begin{cases}
(\frak n;\frak m)&\text{if}\ \frak n>\frak m,\cr
(\frak m;\frak n)&\text{otherwise}.
\end{cases}
\]
We have
\[
\begin{split}
0\le\,&1- \frac{P_k(\frak m;\frak n)}{N_k(\frak m)N_k(\frak n)}\cr
=\,&1- \frac{P_k(\frak m';\frak n')}{N_k(\frak m')N_k(\frak n')}\cr
=\,&\frac{|\{f,g)\in \mathcal N_k(\frak m')\times\mathcal  N_k(\frak n'):\text{gcd}(f,g)\ne 1\}|}{N_k(\frak m')N_k(\frak n')}\cr
\le\,&\frac{|(\mathcal N_k(\frak m')\setminus \mathcal I_k(\frak m'))\times \mathcal N_k(\frak n')|+|\{(f,f)\in \mathcal N_k(\frak m')\times \mathcal N_k(\frak n'):
f\in \mathcal I_k(\frak m')\}|}{N_k(\frak m')N_k(\frak n')}\cr
\le\,&1-\frac{I_k(\frak m')}{N_k(\frak m')}+\frac  1{N_k(\frak m')}\to 0\qquad\text{as $m_{k-1},m_k\to\infty$}.
\end{split}
\]
\end{proof}

We also have the following result which corresponds to Theorem~\ref{T4.2} of the total degree case.

\begin{Thm}\label{T5.5}
Let $k\ge 2$, $\frak m=(m_1,\dots,m_k),\ \frak n=(n_1,\dots,n_k),\ \frak t=(t_1,\dots,t_k)\in\Bbb N^k$ such that $\frak t\le \frak m,\frak n$ and
$\max\{m_i,n_i\}>2t_i+1$ for all $1\le i\le k$. Then
\[
P_k(\frak m;\frak n)=\sum_{\frak 0\le\frak d\le \frak t}N_k(\frak m-\frak d)N_k(\frak n-\frak d)A_k(\frak d)+O\Bigl(q^{\max_{1\le j\le k}(\frac{m_j-t_j}{m_j+1}
\prod_{i=1}^k(m_i+1)+ \frac{n_j-t_j}{n_j+1}\prod_{i=1}^k(n_i+1))}\Bigr).
\]
\end{Thm}

\noindent{\bf Note.} The asymptotic formula in Theorem~\ref{T5.5} is intersting only when \break
$\prod_{i\ne j}(m_i+1)+\prod_{i\ne j}(n_i+1)\gg m_j+n_j$ for all $1\le j\le k$ since otherwise the $O$-term is comparable to $N_k(\frak m)N_k(\frak n)$. In particular, Theorem~\ref{T5.5} does not imply Cororllary 12 of \cite{Coh3}.

\begin{Lem}\label{L5.6}
Let $k\ge 2$,  $\frak m=(m_1,\dots,m_k),\ \frak n=(n_1,\dots,n_k)\in\Bbb N^k$ and $t_1\in\Bbb N$ such that $\max\{m_1,n_1\}>2t_1+1$. Let $\frak d=(d_1,\dots,d_k)\in\Bbb N^k$ such that $\frak d\le \frak m,\frak n$ and $d_1\ge t_1+1$. Then
\[
\begin{split}
&\prod_{i=1}^k(m_i-d_i+1)+\prod_{i=1}^k(n_i-d_i+1)+\prod_{i=1}^k(d_i+1)\cr
\le\,&(m_1-t_1)\prod_{i=2}^k(m_i+1)+(n_1-t_1)\prod_{i=2}^k(n_i+1)-(d_1-(t_1+1))-d_2-\cdots-d_k+t_1+2.
\end{split}
\]
\end{Lem}

\begin{proof}
Let $(\tau_1,\dots,\tau_k)=(t_1+1,0,\dots,0)$.
Let $x_i=m_i-d_i+1$, $y_i=n_i-d_i+1$, $z_i=d_i-\tau_i$. Then $x_i\ge 1$, $y_i\ge 1$, $z_i\ge 0$. Without loss of generality, assume $m_1>2t_1+1$. We have

\[
\begin{split}
&\prod_{i=1}^k(m_i-d_i+1)+\prod_{i=1}^k(n_i-d_i+1)+\prod_{i=1}^k(d_i+1)\cr
&-(m_1-t_1)\prod_{i=2}^k(m_i+1)-(n_1-t_1)\prod_{i=2}^k(n_i+1)\cr
=\,&\prod_{i=1}^k x_i+\prod_{i=1}^k y_i+\prod_{i=1}^k(z_i+\tau_i+1)-\prod_{i=1}^k(x_i+z_i)-\prod_{i=1}^k(y_i+z_i)\cr
=\,&\prod_{i=1}^k x_i+\prod_{i=1}^k y_i+\sum_{I\subset\{1,\dots,k\}}\Bigl(\prod_{i\in I}z_i\Bigr)
\Bigl[\prod_{i\in I'}(\tau_i+1)-\prod_{i\in I'}x_i-\prod_{i\in I'}y_i\Bigr]\qquad (I'=\{1,\dots,k\}\setminus I)\cr
=\,&\tau_1+1+\sum_{\emptyset\ne I\subset\{1,\dots,k\}}\Bigl(\prod_{i\in I}z_i\Bigr)\Bigl[\prod_{i\in I'}(\tau_i+1)-\prod_{i\in I'}x_i-\prod_{i\in I'}y_i\Bigr]\cr
=\,&t_1+2+\sum_{1\in I\subset\{1,\dots,k\}}\Bigl(\prod_{i\in I}z_i\Bigr)\Bigl[1-\prod_{i\in I'}x_i-\prod_{i\in I'}y_i\Bigr]\cr
&+\sum_{\emptyset\ne I\subset\{2,\dots,k\}}\Bigl(\prod_{i\in I}z_i\Bigr)\Bigl[t_1+2-\prod_{i\in I'}x_i-\prod_{i\in I'}y_i\Bigr]\cr
=\,&t_1+2+z_1\Bigl[1-\prod_{i=2}^k x_i-\prod_{i=2}^ky_i\Bigr]+
\sum_{\emptyset\ne I\subset\{2,\dots,k\}}z_1\Bigl(\prod_{i\in I}z_i\Bigr)\Bigl[1-\prod_{i\in I''}x_i-\prod_{i\in I''}y_i\Bigr]\cr
&+\sum_{\emptyset\ne I\subset\{2,\dots,k\}}\Bigl(\prod_{i\in I}z_i\Bigr)\Bigl[t_1+2-x_1\prod_{i\in I''}x_i-y_1\prod_{i\in I''}y_i\Bigr]
\quad(I''=\{2,\dots,k\}\setminus I)\cr
=\,&t_1+2+z_1\Bigl[1-\prod_{i=2}^k x_i-\prod_{i=2}^ky_i\Bigr]\cr
&+\sum_{\emptyset\ne I\subset\{2,\dots,k\}}\Bigl(\prod_{i\in I}z_i\Bigr)\Bigl[t_1+2+z_1-(z_1+x_1)\prod_{i\in I''}x_i-(z_1+y_1)\prod_{i\in I''}y_i\Bigr]
\cr
\le\,&t_1+2-z_1 +\sum_{\emptyset\ne I\subset\{2,\dots,k\}}\Bigl(\prod_{i\in I}z_i\Bigr)\bigl[t_1+2+z_1-(z_1+x_1)-(z_1+y_1)\bigr]\cr
\le\,&t_1+2-z_1 +\sum_{\emptyset\ne I\subset\{2,\dots,k\}}\Bigl(\prod_{i\in I}z_i\Bigr)\bigl[t_1+1-(z_1+x_1)\bigr]\cr
=\,&t_1+2-z_1 -\sum_{\emptyset\ne I\subset\{2,\dots,k\}}\Bigl(\prod_{i\in I}z_i\Bigr)(m_1-2t_1-1)\cr
\le\,&t_1+2-z_1 -\sum_{\emptyset\ne I\subset\{2,\dots,k\}}\prod_{i\in I}z_i\cr
\le\,& t_1+2-z_1-z_2-\cdots-z_k\cr
=\,&t_1+2-(d_1-(t_1+1))-d_2-\cdots-d_k.
\end{split}
\]
\end{proof}

\begin{proof}[Proof of Theorem~\ref{T5.5}]
We have
\[
P_k(\frak m;\frak n)=\Bigl(\sum_{\frak 0\le \frak d\le \frak t}+\sum_{\frak d\not\le\frak t}\Bigr)N_k(\frak m-\frak d)N_k(\frak n-\frak d)A_k(\frak d).
\]
So it suffices to show that
\[
\sum_{\frak d\not\le\frak t}N_k(\frak m-\frak d)N_k(\frak n-\frak d)A_k(\frak d)=O\Bigl(q^{\max_{1\le j\le k}(\frac{m_j-t_j}{m_j+1}
\prod_{i=1}^k(m_i+1)+ \frac{n_j-t_j}{n_j+1}\prod_{i=1}^k(n_i+1))}\Bigr).
\]
Clearly,
\[
\Bigl|\sum_{\frak d\not\le\frak t}N_k(\frak m-\frak d)N_k(\frak n-\frak d)A_k(\frak d)\Bigr|\le
\sum_{j=1}^k\,\sum_{(0,\dots,t_j+1,\dots,0)\le\frak d\le \frak m,\frak n}N_k(\frak m-\frak d)N_k(\frak n-\frak d)N_k(\frak d).
\]
So it suffices to show that
\[
\sum_{(0,\dots,t_j+1,\dots,0)\le\frak d\le \frak m,\frak n}N_k(\frak m-\frak d)N_k(\frak n-\frak d)N_k(\frak d)=O\Bigl(q^{\frac{m_j-t_j}{m_j+1}
\prod_{i=1}^k(m_i+1)+ \frac{n_j-t_j}{n_j+1}\prod_{i=1}^k(n_i+1)}\Bigr).
\]
Without loss of generality let $j=1$. We have
\[
\begin{split}
&\sum_{(t_1+1,0,\dots,0)\le\frak d\le \frak m,\frak n}N_k(\frak m-\frak d)N_k(\frak n-\frak d)N_k(\frak d)\cr
\le\,& \sum_{(t_1+1,0,\dots,0)\le\frak d\le \frak m,\frak n} q^{\prod_{i=1}^k(m_i-d_i+1)+\prod_{i=1}^k(n_i-d_i+1)+\prod_{i=1}^k(d_i+1)}\cr
\le\,& \sum_{(t_1+1,0,\dots,0)\le\frak d\le \frak m,\frak n} q^{(m_1-t_1)\prod_{i=2}^k(m_i+1)+(n_1-t_1)\prod_{i=2}^k(n_i+1)-(d_1-(t_1+1))-d_2-\cdots-d_k+t_1+2}\cr
&\kern10cm \text{(by Lemma~\ref{L5.6})}\cr
=\,&q^{(m_1-t_1)\prod_{i=2}^k(m_i+1)+(n_1-t_1)\prod_{i=2}^k(n_i+1)}q^{t_1+2}\sum_{(t_1+1,0,\dots,0)\le\frak d\le \frak m,\frak n}q^{-(d_1-(t_1+1))-d_2-\cdots-d_k}\cr
=\,&O(q^{(m_1-t_1)\prod_{i=2}^k(m_i+1)+(n_1-t_1)\prod_{i=2}^k(n_i+1)}).
\end{split}
\]
\end{proof}

The next theorem is a variation of Theorem~\ref{T5.5}; it implies Corollary 12 of \cite{Coh3}.

\vskip 1cm

\begin{Thm}\label{T5.7}
Let $k\ge 2$, $\frak m=(m_1,\dots,m_k), \frak n=(n_1,\dots,n_k)\in\Bbb N^k$ such that $m_k>0$, $n_k>0$ and $\max\{m_i,n_i\}>1$ for all $1\le i\le k-1$. Then
\[
\begin{split}
P_k(\frak m;\frak n)=\,&N_k(\frak m)N_k(\frak n)-qN_k(m_1,\dots,m_{k-1},m_k-1)N_k(n_1,\dots,n_{k-1},n_k-1)\cr
&+O\Bigl(q^{\max_{1\le j\le k-1}(\frac{m_j}{m_j+1}\prod_{i=1}^k(m_i+1)+\frac{n_j}{n_j+1}\prod_{i=1}^k(n_i+1))}\Bigr).
\end{split}
\]
\end{Thm}

\begin{proof}
We have
\[
P_k(\frak m;\frak n)=\Bigl(\sum_{\substack{\frak d=(0,\dots,0,d_k)\cr 0\le d_k\le m_k,n_k}}+
\sum_{\substack{\frak d=(d_1,\dots,d_k)\le \frak m,\frak n\cr (d_1,\dots,d_{k-1})\ne\frak 0}}\Bigr)N_k(\frak m-\frak d)N_k(\frak n-\frak d)A_k(\frak d).
\]
By \eqref{A1d},
\[
\begin{split}
&\sum_{\substack{\frak d=(0,\dots,0,d_k)\cr 0\le d_k\le m_k,n_k}}N_k(\frak m-\frak d)N_k(\frak n-\frak d)A_k(\frak d)\cr
=\,&
N_k(\frak m)N_k(\frak n)-qN_k(m_1,\dots,m_{k-1},m_k-1)N_k(n_1,\dots,n_{k-1},n_k-1).
\end{split}
\]
By the proof of Theorem~\ref{T5.5}, we also have
\[
\Bigl|\sum_{\substack{\frak d=(d_1,\dots,d_k)\le \frak m,\frak n\cr (d_1,\dots,d_{k-1})\ne\frak 0}}N_k(\frak m-\frak d)N_k(\frak n-\frak d)A_k(\frak d)\Bigr|=
O\Bigl(q^{\max_{1\le j\le k-1}(\frac{m_j}{m_j+1}\prod_{i=1}^k(m_i+1)+\frac{n_j}{n_j+1}\prod_{i=1}^k(n_i+1))}\Bigr).
\]
Therefore the conclusion follows.
\end{proof}


\section*{Appendix}

\appendix


\section{Proof of Lemma~\ref{L0.3}}

We need two additional lemmas for the proof of Lemma~\ref{L0.3}

\begin{Lem}\label{LB.1}
Let $(m_1^{(j)},\dots,m_k^{(j)})\in \Bbb N^k$, $1\le j\le l$. Then
\[
\prod_i\Bigl(\sum_j m_i^{(j)}+1\Bigr)-\sum_j\prod_i(m_i^{(j)}+1)\ge-(l-1).
\]
\end{Lem}

\begin{proof}
Using induction on $l$ we only have to prove the case $l=2$, i.e.,
\[
\prod_i(m_i+n_i+1)-\prod_i(m_i+1)-\prod_i(n_i+1)\ge -1,
\]
where $m_i,n_i\in\Bbb N$.
If all nonzero entries of
\[
\left[
\begin{matrix}
m_1&\cdots& m_k\cr
n_1&\cdots& n_k
\end{matrix}\right]
\]
appear in a single column, say $\left[\begin{smallmatrix}
m_2&\cdots& m_k\cr
n_2&\cdots& n_k
\end{smallmatrix}\right]=0$, the conclusion is obviously true. So assume $m_1>0$ and $n_2>0$. Then
\[
(m_1+n_1+1)(m_2+n_2+1)-(m_1+1)(m_2+1)-(n_1+1)(n_2+1)=m_1n_2+m_2n_1-1\ge 0.
\]
So
\[
\begin{split}
&\prod_i(m_i+n_i+1)-\prod_i(m_i+1)-\prod_i(n_i+1)\cr
\ge\,&\bigl[(m_1+n_1+1)(m_2+n_2+1)-(m_1+1)(m_2+1)-(n_1+1)(n_2+1)\bigr]\prod_{i=3}^k(m_i+n_i+1)\cr
\ge\,&0.
\end{split}
\]
\end{proof} 

\begin{Lem}\label{LB.2}
Let $k\ge 2$ and let $(m_1,\dots,m_k),(n_1,\dots,n_k)\in\Bbb N^k$ such that $m_i+n_i>0$ for all $1\le i\le k$. For each $a\in\Bbb R$, write $\overline a=\max\{a,0\}$. Then
\[
\begin{split}
&\prod_i(m_i+n_i+1)-\prod_i(m_i+1)-\prod_i(n_i+1)\cr
&
\begin{cases}
=-1\kern4.7cm\text{if $(m_1,\dots,m_k)=0$ or $(n_1,\dots,n_k)=0$},\cr \cr
\ge \min\{\overline{\frac 12(m_1-1)},\overline{\frac 12(n_1-1)}\}+\min\{\overline{\frac 12(m_2-1)},\overline{\frac 12(n_2-1)}\}+m_3+\cdots+m_k\cr
\kern5.4cm\text{if $(m_1,\dots,m_k)\ne0$ and $(n_1,\dots,n_k)\ne0$}.
\end{cases}
\end{split}
\]
\end{Lem}

\begin{proof}
We only have to prove the claim when $(m_1,\dots, m_k)\ne 0$ and $(n_1,\dots,n_k)\ne 0$. First assume $k=2$.  We have
\[
(m_1+n_1+1)(n_2+n_2+1)-(m_1+1)(m_2+1)-(n_1+1)(n_2+1)=m_1n_2+m_2n_1-1.
\]
Assume $m_1\ge m_2$.
If $n_1=0$ or $m_2=0$, then $m_1,n_2>0$. So
\[
m_1n_2+m_2n_1-1\ge m_1-1\ge {\textstyle\frac 12(m_1-1)}+\overline{\textstyle\frac 12(m_2-1)}.
\]
If $n_1>0$ and $m_2>0$, then
\[
m_1n_2+m_2n_1-1\ge m_2n_1-1\ge (n_1-1)+(m_2-1).
\]

Now assume $k\ge3$. 

{\bf Case 1.} Assume that there exist $1\le a<b\le k$ such that $m_an_b+m_bn_a-1>0$, say $m_1n_2+m_2n_1-1>0$. Then
\[
\begin{split}
&\prod_i(m_i+n_i+1)-\prod_i(m_i+1)-\prod_i(n_i+1)\cr
\ge\,&(m_1n_2+m_2n_1-1)\prod_{i=3}^k(m_i+n_i+1)\cr
\ge\,&(m_1n_2+m_2n_1-1)+\prod_{i=3}^k(m_i+n_i+1)-1\cr
\ge\,&m_1n_2+m_2n_1-1+m_3+\cdots+m_k\cr
\ge\,&\min\{\overline{\textstyle\frac 12(m_1-1)},\overline{\textstyle\frac 12(n_1-1)}\}+\min\{\overline{\textstyle\frac 12(m_2-1)},\overline{\textstyle\frac 12(n_2-1)}\}+m_3+\cdots+m_k
\end{split}
\]

{\bf Case 2.} Assume that there do not exist $1\le a<b\le k$ such that $m_an_b+m_bn_a-1>0$. Then
\[
\left[
\begin{matrix} m_1&\cdots&m_k\cr n_1&\cdots&n_k\end{matrix}\right]=
\left[\phantom{\begin{matrix} 0\cr 0\end{matrix}}\right. \kern-2mm
\underbrace{\begin{matrix} 1&\cdots&1\cr 0&\cdots&0\end{matrix}}_{\alpha>0}\ \underbrace{\begin{matrix}0&\cdots&0\cr 1&\cdots&1\end{matrix}}_{\beta>0}
\kern -2mm\left.\phantom{\begin{matrix} 0\cr 0\end{matrix}}\right]
\quad \text{or}\quad
\left[\phantom{\begin{matrix} 0\cr 0\end{matrix}}\right. \kern-2mm
\underbrace{\begin{matrix} 1&\cdots&1\cr 0&\cdots&0\end{matrix}}_{\alpha'}\ \begin{matrix} 1\cr 1\end{matrix}\ 
\underbrace{\begin{matrix}0&\cdots&0\cr 1&\cdots&1\end{matrix}}_{\beta'}
\kern -2mm\left.\phantom{\begin{matrix} 0\cr 0\end{matrix}}\right].
\]
In the first case,
\[
\begin{split}
&\prod_{i=1}^k(m_i+n_i+1)-\prod_{i=1}^k(m_i+1)-\prod_{i=1}^k(n_i+1)\cr
=\,&2^k-2^\alpha-2^\beta\le 2^k-2^{k-1}-2=2^{k-1}-2\ge k-1\ge \alpha=m_1+\cdots+m_k.
\end{split}
\] 
In the second case,
\[
\begin{split}
&\prod_{i=1}^k(m_i+n_i+1)-\prod_{i=1}^k(m_i+1)-\prod_{i=1}^k(n_i+1)\cr
=\,&3\cdot 2^{k-1}-2^{\alpha'+1}-2^{\beta'+1}\le 3\cdot2^{k-1}-2^k-2=2^{k-1}-2\cr
\ge\,& k-1\ge \overline{\textstyle\frac 12(m_1-1)}+m_2+\cdots+m_k.
\end{split}
\] 
\end{proof}

\begin{proof}[Proof of Lemma~\ref{L0.3}]
As mentioned before, the sum in Lemma~\ref{L0.3} is the cardinality of 
\[
\mathcal F=\{f\in \mathcal N_k(\frak m): \text{all irreducible factors of $f$ have ${\deg}_{x_1}\le m_1-t_1$}\}.
\]
Let
\[
\mathcal F_1=\{:f\in \mathcal F:f=f_1f_2,\ t_1<{\deg}_{x_1}f_1\le m_1-t_1\}.
\]
Then it suffices to show that 
\begin{gather}
\label{26}
|\mathcal F_1|=O(q^{(m_1-t_1+1)(m_2+1)\cdots(m_k+1)}),\\ 
\label{27}
|\mathcal F\setminus\mathcal F_1|=O(q^{(m_1-t_1+1)(m_2+1)\cdots(m_k+1)}). 
\end{gather}

We claim that 
\begin{equation}\label{28}
\mathcal F\setminus\mathcal F_1\subset\{f\in\mathcal F: f=f_1f_2f_3,\ {\deg}_{x_1}f_i\le t_1,\ i=1,2,3\}.
\end{equation}
Let $f\in\mathcal F\setminus\mathcal F_1$. Then all irreducible factors of $f$ have $\deg_{x_1}\le t_1$. Write $f=f_1f_2f_3$ where $\deg_{x_1}f_1\le\deg_{x_1}f_2\le\deg_{x_1}f_3$ such that $\deg_{x_1}f_3$ is as small as possible. Let $d_i=\deg_{x_1}f_i$. If $d_3\le t_1$, we are done. So assume $d_3>t_1$. We must have $d_3>m_1-t_1$ since otherwise
$f\in\mathcal F_1$. Then $f_3=f_3'f_3''$ where $\deg_{x_1}f_3'<d_3$ and $\deg_{x_3}f_3''<d_3$. We may assume $\deg_{x_1}f_3'\le\frac {d_3}2$. Then
\[
f=(f_1f_3')\cdot f_2\cdot f_3''
\]
where
\begin{gather*}
{\deg}_{x_1}(f_1f_3')=d_1+\frac{d_3}2\le \frac{d_1+d_2}2+\frac {d_3}2=\frac {m_1}2\le m_1-t_1<d_3,\\
{\deg}_{x_1}f_2\le m_1-d_3<t_1<d_3,\\
{\deg}_{x_1}f_3''<d_3.
\end{gather*}
This contradicts the minimality of $d_3$. So \eqref{28} is proved.

\medskip

{\bf Case 1.} Assume that $m_1>2t_1$.

By \eqref{28} we have
\[
\begin{split}
&|\mathcal F\setminus\mathcal F_1|\cr
\le\,&\sum_{\sum_{j=1}^3(m_2^{(j)},\dots,m_k^{(j)})=(m_2,\dots,m_k)}\prod_{j=1}^3 N_k(t_1,m_2^{(j)},\dots,m_k^{(j)})\cr
\le\,&\sum_{\sum_{j=1}^3(m_2^{(j)},\dots,m_k^{(j)})=(m_2,\dots,m_k)}q^{\sum_{j=1}^3(t_1+1)\prod_{i=2}^k(m_i^{(j)}+1)}\cr
\le\,&\sum_{\sum_{j=1}^3(m_2^{(j)},\dots,m_k^{(j)})=(m_2,\dots,m_k)}q^{(t_1+1)[\prod_{i=2}^k(m_i+1)+2]} \qquad\text{(by Lemma~\ref{LB.1})}\cr
=\,&q^{2(t_1+1)}\sum_{\sum_{j=1}^3(m_2^{(j)},\dots,m_k^{(j)})=(m_2,\dots,m_k)}q^{(t_1+1)\prod_{i=2}^k(m_i+1)}\cr
\le\,&q^{2(t_1+1)}\sum_{\sum_{j=1}^3(m_2^{(j)},\dots,m_k^{(j)})=(m_2,\dots,m_k)}q^{(m_1-t_1+1)\prod_{i=2}^k(m_i+1)-\prod_{i=2}^k(m_i+1)}
\qquad (\because m_1-t_1>t_1)\cr
\le\,&q^{2(t_1+1)}q^{(m_1-t_1+1)\prod_{i=2}^k(m_i+1)}\sum_{\sum_{j=1}^3(m_2^{(j)},\dots,m_k^{(j)})=(m_2,\dots,m_k)}q^{-m_2-\cdots-m_k}\cr
\le\,&q^{2(t_1+1)}q^{(m_1-t_1+1)\prod_{i=2}^k(m_i+1)}\sum_{m_i^{(i)},\, 2\le i\le k,\,j=1,2,3}q^{-\sum_{i,j}m_i^{(i)}}\cr
=\,&O(q^{(m_1-t_1+1)\prod_{i=2}^k(m_i+1)}).
\end{split}
\]

Now we prove \eqref{26}. 
Write 
\[
\mathcal F_1=\mathcal F_1'\cup\mathcal F_1''
\]
where
\begin{gather*}
\mathcal F_1'=\{f\in\mathcal F_1: f=f_1f_2,\ t_1<{\deg}_{x_1}f_1\le m_1-t_1,\ ({\deg}_{x_2}f_i,\dots,{\deg}_{x_k}f_i)\ne0 \ \text{for $i=1$ and $2$}\},\\
\mathcal F_1''=\{f\in\mathcal F_1: f=f_1f_2,\ t_1<{\deg}_{x_1}f_1\le m_1-t_1,\ ({\deg}_{x_2}f_i,\dots,{\deg}_{x_k}f_i)=0 \ \text{for $i=1$ or $2$}\}.
\end{gather*}
We prove in turn that both $|\mathcal F_1''|$ and $|\mathcal F_1'|$ are $O(q^{(m_1-t_1+1)\prod_{i+2}^k(m_i+1)})$.

We have 
\[
\begin{split}
&|\mathcal F_1''|\cr
\le\,&\sum_{t_1<d\le m_1-t_1}\bigl(N_k(d,0,\dots,0)N_k(m_1-d,m_2,\dots,m_k)+N_k(d,m_2,\dots,m_k)N_k(m_1-d,0,\dots,0)\bigr)\cr
\le\,&\sum_{t_1<d\le m_1-t_1}(q^{(d+1)+(m_1-d+1)\prod_{i=2}^k(m_i+1)}+q^{(d+1)\prod_{i=2}^k(m_i+1)+(m_1-d+1)}),
\end{split}
\]
where
\begin{equation}\label{29}
\begin{split}
&\sum_{t_1<d\le m_1-t_1}q^{(d+1)+(m_1-d+1)\prod_{i=2}^k(m_i+1)}\cr
=\,&q^{(m_1-t_1+1)\prod_{i=2}^k(m_i+1)}\sum_{t_1<d\le m_1-t_1}q^{(d+1)-(d-t_1)\prod_{i=2}^k(m_i+1)}\cr
\le\,&q^{(m_1-t_1+1)\prod_{i=2}^k(m_i+1)}\sum_{t_1<d\le m_1-t_1}q^{(d+1)-2(d-t_1)}\cr
=\,&q^{(m_1-t_1+1)\prod_{i=2}^k(m_i+1)}q^{t_1+1}\sum_{t_1<d\le m_1-t_1}q^{-(d-t_1)}\cr
=\,&O(q^{(m_1-t_1+1)\prod_{i+2}^k(m_i+1)})
\end{split}
\end{equation}
and
\[
\begin{split}
&\sum_{t_1<d\le m_1-t_1}q^{(d+1)\prod_{i=2}^k(m_i+1)+(m_1-d+1)}\cr
=\,&q^{(m_1-t_1+1)\prod_{i=2}^k(m_i+1)+(t_1+1)}+\sum_{t_1<d< m_1-t_1}q^{(d+1)\prod_{i=2}^k(m_i+1)+(m_1-d+1)}\cr
=\,&q^{(m_1-t_1+1)\prod_{i=2}^k(m_i+1)}q^{(t_1+1)}+\sum_{t_1<d'< m_1-t_1}q^{(d'+1)+(m_1-d'+1)\prod_{i=2}^k(m_i+1)}\qquad (d'=m_1-d)\cr
=\,&O(q^{(m_1-t_1+1)\prod_{i=2}^k(m_i+1)})\qquad\text{(by \eqref{29})}.
\end{split}
\]

It remains to show that $|\mathcal F_1'|=O(q^{(m_1-t_1+1)\prod_{i+2}^k(m_i+1)})$. We first assume $k\ge 3$. 
Write $\frak m^{(j)}=(m_2^{(j)},\dots,m_k^{(j)})\in\Bbb N^{k-1}$, $j=1,2$.
We have
\[
\begin{split}
&|\mathcal F_1'|\cr
\le\,&\sum_{t_1<d\le m_1-t_1}\sum_{\substack{\frak m^{(1)},\frak m^{(2)}\ne\frak 0\cr \frak m^{(1)}+\frak m^{(2)}=(m_2,\dots, m_k)}}
N_k(d,\frak m^{(1)})N_k(m_1-d,\frak m^{(2)})\cr
\le\,&\sum_{t_1<d\le m_1-t_1}\sum_{\substack{\frak m^{(1)},\frak m^{(2)}\ne\frak 0\cr \frak m^{(1)}+\frak m^{(2)}=(m_2,\dots, m_k)}}q^{(d+1)\prod_{i=2}^k(m_i^{(1)}+1)+(m_1-d+1)\prod_{i=2}^k(m_i^{(2)}+1)}\cr
\le\,&\sum_{t_1<d\le m_1-t_1}\sum_{\substack{\frak m^{(1)},\frak m^{(2)}\ne\frak 0\cr \frak m^{(1)}+\frak m^{(2)}=(m_2,\dots, m_k)}}q^{(m_1-t_1+1)[\prod_{i=2}^k(m_i^{(1)}+1)+\prod_{i=2}^k(m_i^{(2)}+1)]-(d-t_1)\prod_{i=2}^k(m_i^{(2)}+1)}\cr
\le\,&\sum_{t_1<d\le m_1-t_1}\sum_{\substack{\frak m^{(1)},\frak m^{(2)}\ne\frak 0\cr \frak m^{(1)}+\frak m^{(2)}=(m_2,\dots, m_k)}}q^{(m_1-t_1+1)\prod_{i=2}^k(m_i+1)}q^{-(d-t_1)\prod_{i=2}^k(m_i^{(2)}+1)}
\qquad\text{(by Lemma~\ref{LB.2})}\cr
\le\,&q^{(m_1-t_1+1)\prod_{i=2}^k(m_i+1)}\;\sum_{t_1<d\le m_1-t_1}\sum_{m_i^{(2)},\; 2\le i\le k}q^{-(d-t_1)-m_2^{(2)}-\cdots-m_k^{(2)}}\cr
=\,&O(q^{(m_1-t_1+1)\prod_{i+2}^k(m_i+1)}).
\end{split}
\]

Now assume $k=2$. We have
\[
\begin{split}
&|\mathcal F_1'|\cr
\le\,&\sum_{t_1<d\le m_1-t_1}\;\sum_{a+b=m_2}q^{(d+1)(a+1)+(m_1-d+1)(b+1)}\cr
=\,&q^{(m_1-t_1+1)(m_2+1)}q^{t_1+1}\sum_{t_1<d\le m_1-t_1}\;\sum_{a+b=m_2}q^{-(d-t_1)b-(m_1-t_1-d)a}\cr
\le\,&q^{(m_1-t_1+1)(m_2+1)}q^{t_1+1}\Bigl(\sum_{t_1<d\le m_1-t_1}\;\sum_{b>0}q^{-(d-t_1)b}+\sum_{t_1<d\le m_1-t_1}q^{-(m_1-t_1-d)m_2}\Bigr)\cr
\le\,&q^{(m_1-t_1+1)(m_2+1)}q^{t_1+1}\Bigl(\sum_{t_1<d\le m_1-t_1}\;\sum_{b>0}q^{-(d-t_1)-b+1}+\sum_{t_1<d\le m_1-t_1}q^{-(m_1-t_1-d)}\Bigr)\cr
=\,&O(q^{(m_1-t_1+1)\prod_{i+2}^k(m_i+1)}).
\end{split}
\]
This completes the proof in case 1.

\medskip

{\bf Case 2.} Assume that $m_1=2t_1$ and $k\ge 3$. In this case 
$\mathcal F_1=\emptyset$, so it suffices to prove \eqref{27}.
By \eqref{28} we have
\[
\mathcal F\setminus\mathcal F_1=\mathcal A\cup\mathcal B
\]
where
\begin{gather*}
\mathcal A=\{f\in \mathcal F:f=f_1f_2f_3,\ {\deg}_{x_1}f_i<t_1,\ i=1,2,3\},\\
\mathcal B=\{f\in\mathcal F:f=f_1f_2,\ {\deg}_{x_1}f_1={\deg}_{x_1}f_2=t_1\}.
\end{gather*}
The proof that $|\mathcal A|=O(q^{(t_1+1)\prod_{i=2}^k(m_i+1)})$ is the same 
as the proof that $|\mathcal F\setminus\mathcal F_1|=O(q^{(t_1+1)\prod_{i=2}^k(m_i+1)})$
in case 1. As for $|\mathcal B|$, we have
\[
\begin{split}
&|\mathcal B|\cr
\le\,&\sum_{(a_2,\dots,a_k)+(b_2,\dots,b_k)=(m_2,\dots,m_k)}N_k(t_1,a_2,\dots,a_k)N_k(t_1,b_2,\dots,b_k)\cr
\le\,&2q^{(t_1+1)\prod_{i=2}^k(m_i+1)}+\sum_{\substack{(a_2,\dots,a_k)+(b_2,\dots,b_k)=(m_2,\dots,m_k)\cr
(a_2,\dots,a_k),(b_2,\dots,b_k)\ne 0}}q^{(t_1+1)[\prod_{i=2}^k(a_i+1)+\prod_{i=2}^k(b_i+1)]}\cr
\le\,&2q^{(t_1+1)\prod_{i=2}^k(m_i+1)}+q^{(t_1+1)\prod_{i=2}^k(m_i+1)}\cr
&\cdot\sum_{\substack{(a_2,\dots,a_k)+(b_2,\dots,b_k)=(m_2,\dots,m_k)\cr
(a_2,\dots,a_k),(b_2,\dots,b_k)\ne 0}}q^{-\min\{\overline{\frac 12(a_2-1)},\overline{\frac 12(b_2-1)}\}-\min\{\overline{\frac 12(a_3-1)},\overline{\frac 12(b_3-1)}\}-a_4-\cdots-a_k}\cr
&\kern9cm \text{(by Lemma~\ref{LB.2})}\cr
=\,&O(q^{(t_1+1)\prod_{i=2}^k(m_i+1)}).
\end{split}
\]
This completes the proof of the lemma.
\end{proof}


\section{Tables}


Table~\ref{Tb4} contains the values of $P_2(m_1,m_2;n_1,n_2)$ with $q=2$ and $(m_1,m_2),(n_1,n_2)\le (4,4)$. To present the data efficiently, we observe that 
$P_2(m_1,m_2;n_1,n_2)$ is invariant under row and column permutations of $\left[\begin{smallmatrix} m_1&m_2\cr n_1&n_2\end{smallmatrix}\right]$. Under row and column permutations, matrices  $\left[\begin{smallmatrix} m_1&m_2\cr n_1&n_2\end{smallmatrix}\right]$, $0\le m_1,m_2,n_1,n_2\le 4$,
are represented by
\[
\begin{split}
\left[\begin{matrix} a&b\cr c&d\end{matrix}\right]\ :\ & 0\le a<b,c,d\le 4,\cr
\left[\begin{matrix} a&a\cr b&c\end{matrix}\right]\ :\ & 0\le a<b\le c\le 4,\cr
\left[\begin{matrix} a&b\cr a&c\end{matrix}\right]\ :\ & 0\le a<b\le c\le 4,\cr
\left[\begin{matrix} a&b\cr c&a\end{matrix}\right]\ :\ & 0\le a<b\le c\le 4,\cr
\left[\begin{matrix} a&a\cr a&b\end{matrix}\right]\ :\ & 0\le a\le b\le 4.
\end{split}
\]

\begin{table}
\caption{$I_2(m),\ q=2,3,4,5,\ 0\le m\le 10$}\label{Tb1}
\[
\begin{tabular}{c|l|l}
\hline
$m\backslash q$ & \hfil 2&  \hfil 3\\ \hline
0&0&0 \\
1 & 6 & 12 \\
2 & 35 & 273 \\
3 & 694 & 25520 \\
4 & 26089 & 6778629 \\
5 & 1862994 & 5132148528 \\
6 & 253247715 & 11368775698280 \\
7 & 66799608630 & 74897449398451680 \\
8 & 34698378752226 & 1476178370884382958936 \\
9 & 35781375988234520 & 87205387550224830516286800 \\
10 & 73534241823793715433 & 15450442981642705273095610563240 \\
\hline
\end{tabular}
\]

\medskip
\[
\begin{tabular}{c|l}
\hline
$m\backslash q$ & \hfil 4\\ \hline
0&0 \\
1 & 20 \\
2 & 1134 \\
3 & 323940 \\
4 & 350195076 \\
5 & 1458203653116 \\
6 & 23988036558291750 \\
7 & 1573616297933972778420 \\
8 & 412613600502090075171985440 \\
9 & 432682737835397726783364117773760 \\
10 & 1814830203343733351868975985798075240938 \\
\hline
\end{tabular}
\]

\medskip
\[
\begin{tabular}{c|l}
\hline
$m\backslash q$ & \hfil 5\\ \hline
0&0 \\
1 & 30 \\
2 & 3410 \\
3 & 2330240 \\
4 & 7549603600 \\
5 & 118965950703744 \\
6 & 9309505329218297280 \\
7 & 3637689729211851543816960 \\
8 & 7105314552536912564123328420000 \\
9 & 69388718760088702173445263653542192000 \\
10 & 3388129637939157475672361687005401831354725568 \\
\hline
\end{tabular}
\]
\end{table}
\vskip5mm


\begin{table}
\caption{$P_2(m;n)$ v.s. $N_2(m)N_2(n)$, $q=2$, $m,n\le 5$}\label{Tb2}
\[
\begin{tabular}{c||l|l|l|l|l|l}
\hline
$m\diagdown n$ & \hfil 0&\hfil 1&\hfil 2&\hfil 3&\hfil 4&\hfil 5\\
\hline\hline
0 & 1 & & & & & \\
 & 1 & & & & & \\ \hline
1 & 6& 30 & & & &\\
 & 6& 36 & & & &\\ \hline
2 & 56& 300 & 2900 & & &\\
 & 56& 336 & 3136 & & &\\ \hline
3 & 960& 5424 & 51624 & 901560 & & \\
 & 960& 5760 & 53760 & 921600 & & \\ \hline
4 & 31744& 184704 & 1741984 & 30141936 & 1002049232 & \\
 & 31744& 190464 & 1777664 & 30474240 & 1007681536 & \\ \hline
5 & 2064384& 12195840 & 114443520 & 1970999232 & 65347584672 & 4255612716000\\
 & 2064384& 12386304 & 115605504 & 1981808640 & 65531805696 & 4261681299456\\ \hline
\multicolumn{7}{l}{ }\\
\multicolumn{7}{l}{In each entry, the top number is $P_2(m;n)$; the bottom number is $N_2(m)N_2(n)$.}\\
\end{tabular}
\]
\end{table}


\begin{table}
\caption{$I_2(m_1,m_2)$, $q=2$, $m_1\le m_2\le 5$}\label{Tb3}
\[
\begin{tabular}{|c||l||}
\hline
$m_1,m_2$ & $I_2(m_1,m_2)$\\ \hline\hline
0 , 0  &  0 \\
0 , 1  &  2 \\
0 , 2  &  1 \\
0 , 3  &  2 \\
0 , 4  &  3 \\
0 , 5  &  6 \\
1 , 1  &  6 \\
\hline
\end{tabular}
\begin{tabular}{|c||l||}
\hline
$m_1,m_2$ & $I_2(m_1,m_2)$\\ \hline\hline
1 , 2  &  24 \\
1 , 3  &  96 \\
1 , 4  &  384 \\
1 , 5  &  1536 \\
2 , 2  &  243 \\
2 , 3  &  2256 \\
2 , 4  &  19476 \\
\hline
\end{tabular}
\begin{tabular}{|c||l|}
\hline
$m_1,m_2$ & $I_2(m_1,m_2)$\\ \hline\hline
2 , 5  &  162816 \\
3 , 3  &  43798 \\
3 , 4  &  774240 \\
3 , 5  &  13042176 \\
4 , 4  &  27518145 \\
4 , 5  &  927161664 \\
5 , 5  &  62409885906 \\
\hline
\end{tabular}
\]
\end{table}


\begin{table}
\caption{$P_2(\frak m;\frak n)$, $q=2$, $\frak m, \frak n\le (4,4)$}\label{Tb4}
\noindent $\left[\begin{smallmatrix} \frak m\cr \frak n\end{smallmatrix}\right]=\left[\begin{smallmatrix} a&b\cr c&d\end{smallmatrix}\right]:
0\le a<b,c,d\le 4$
\[
\begin{tabular}{|c||l||}
\hline
$a,b,c,d$ & $P_2(a,b;c,d)$\\ \hline\hline
0 , 1 , 1 , 1  &  16 \\
0 , 1 , 1 , 2  &  68 \\
0 , 1 , 1 , 3  &  280 \\
0 , 1 , 1 , 4  &  1136 \\
0 , 1 , 2 , 1  &  80 \\
0 , 1 , 2 , 2  &  712 \\
0 , 1 , 2 , 3  &  5984 \\
0 , 1 , 2 , 4  &  49024 \\
0 , 1 , 3 , 1  &  352 \\
0 , 1 , 3 , 2  &  6416 \\
0 , 1 , 3 , 3  &  108928 \\
0 , 1 , 3 , 4  &  1793024 \\
0 , 1 , 4 , 1  &  1472 \\
0 , 1 , 4 , 2  &  54304 \\
0 , 1 , 4 , 3  &  1852928 \\
0 , 1 , 4 , 4  &  61136896 \\
0 , 2 , 1 , 1  &  32 \\
0 , 2 , 1 , 2  &  136 \\
0 , 2 , 1 , 3  &  560 \\
0 , 2 , 1 , 4  &  2272 \\
0 , 2 , 2 , 1  &  160 \\
0 , 2 , 2 , 2  &  1424 \\
0 , 2 , 2 , 3  &  11968 \\
0 , 2 , 2 , 4  &  98048 \\
0 , 2 , 3 , 1  &  704 \\
0 , 2 , 3 , 2  &  12832 \\
0 , 2 , 3 , 3  &  217856 \\
0 , 2 , 3 , 4  &  3586048 \\
0 , 2 , 4 , 1  &  2944 \\
0 , 2 , 4 , 2  &  108608 \\
0 , 2 , 4 , 3  &  3705856 \\
0 , 2 , 4 , 4  &  122273792 \\
0 , 3 , 1 , 1  &  64 \\
0 , 3 , 1 , 2  &  272 \\
\hline
\end{tabular}
\begin{tabular}{|c||l||}
\hline
$a,b,c,d$ & $P_2(a,b;c,d)$\\ \hline\hline
0 , 3 , 1 , 3  &  1120 \\
0 , 3 , 1 , 4  &  4544 \\
0 , 3 , 2 , 1  &  320 \\
0 , 3 , 2 , 2  &  2848 \\
0 , 3 , 2 , 3  &  23936 \\
0 , 3 , 2 , 4  &  196096 \\
0 , 3 , 3 , 1  &  1408 \\
0 , 3 , 3 , 2  &  25664 \\
0 , 3 , 3 , 3  &  435712 \\
0 , 3 , 3 , 4  &  7172096 \\
0 , 3 , 4 , 1  &  5888 \\
0 , 3 , 4 , 2  &  217216 \\
0 , 3 , 4 , 3  &  7411712 \\
0 , 3 , 4 , 4  &  244547584 \\
0 , 4 , 1 , 1  &  128 \\
0 , 4 , 1 , 2  &  544 \\
0 , 4 , 1 , 3  &  2240 \\
0 , 4 , 1 , 4  &  9088 \\
0 , 4 , 2 , 1  &  640 \\
0 , 4 , 2 , 2  &  5696 \\
0 , 4 , 2 , 3  &  47872 \\
0 , 4 , 2 , 4  &  392192 \\
0 , 4 , 3 , 1  &  2816 \\
0 , 4 , 3 , 2  &  51328 \\
0 , 4 , 3 , 3  &  871424 \\
0 , 4 , 3 , 4  &  14344192 \\
0 , 4 , 4 , 1  &  11776 \\
0 , 4 , 4 , 2  &  434432 \\
0 , 4 , 4 , 3  &  14823424 \\
0 , 4 , 4 , 4  &  489095168 \\
1 , 2 , 2 , 2  &  16304 \\
1 , 2 , 2 , 3  &  139480 \\
1 , 2 , 2 , 4  &  1152656 \\
1 , 2 , 3 , 2  &  142144 \\
\hline
\end{tabular}
\begin{tabular}{|c||l|}
\hline
$a,b,c,d$ & $P_2(a,b;c,d)$\\ \hline\hline
1 , 2 , 3 , 3  &  2448560 \\
1 , 2 , 3 , 4  &  40593472 \\
1 , 2 , 4 , 2  &  1184768 \\
1 , 2 , 4 , 3  &  40955488 \\
1 , 2 , 4 , 4  &  1360009472 \\
1 , 3 , 2 , 2  &  68896 \\
1 , 3 , 2 , 3  &  585296 \\
1 , 3 , 2 , 4  &  4820800 \\
1 , 3 , 3 , 2  &  601088 \\
1 , 3 , 3 , 3  &  10288288 \\
1 , 3 , 3 , 4  &  170027072 \\
1 , 3 , 4 , 2  &  5012224 \\
1 , 3 , 4 , 3  &  172190528 \\
1 , 3 , 4 , 4  &  5700390784 \\
1 , 4 , 2 , 2  &  282944 \\
1 , 4 , 2 , 3  &  2396320 \\
1 , 4 , 2 , 4  &  19705088 \\
1 , 4 , 3 , 2  &  2469376 \\
1 , 4 , 3 , 3  &  42142016 \\
1 , 4 , 3 , 4  &  695422336 \\
1 , 4 , 4 , 2  &  20595200 \\
1 , 4 , 4 , 3  &  705500800 \\
1 , 4 , 4 , 4  &  23322302720 \\
2 , 3 , 3 , 3  &  192236480 \\
2 , 3 , 3 , 4  &  3180273904 \\
2 , 3 , 4 , 3  &  3193241536 \\
2 , 3 , 4 , 4  &  105804817760 \\
2 , 4 , 3 , 3  &  1586111872 \\
2 , 4 , 3 , 4  &  26194635488 \\
2 , 4 , 4 , 3  &  26352940928 \\
2 , 4 , 4 , 4  &  871678913728 \\
3 , 4 , 4 , 4  &  29920251144512 \\
 & \\
 & \\
\hline
\end{tabular}
\]
\end{table}
\addtocounter{table}{-1}
\begin{table}
\caption{(continued)}
\noindent $\left[\begin{smallmatrix} \frak m\cr \frak n\end{smallmatrix}\right]=\left[\begin{smallmatrix} a&a\cr b&c\end{smallmatrix}\right]:
0\le a<b\le c\le 4$
\[
\begin{tabular}{|c||l||}
\hline
$a,b,c$ & $P_2(a,a;b,c)$\\ \hline\hline
0 , 1 , 1  &  10 \\
0 , 1 , 2  &  44 \\
0 , 1 , 3  &  184 \\
0 , 1 , 4  &  752 \\
0 , 2 , 2  &  400 \\
0 , 2 , 3  &  3392 \\
0 , 2 , 4  &  27904 \\
\hline
\end{tabular}
\begin{tabular}{|c||l||}
\hline
$a,b,c$ & $P_2(a,a;b,c)$\\ \hline\hline
0 , 3 , 3  &  57856 \\
0 , 3 , 4  &  954368 \\
0 , 4 , 4  &  31522816 \\
1 , 2 , 2  &  3628 \\
1 , 2 , 3  &  31496 \\
1 , 2 , 4  &  262096 \\
1 , 3 , 3  &  550624 \\
\hline
\end{tabular}
\begin{tabular}{|c||l|}
\hline
$a,b,c$ & $P_2(a,a;b,c)$\\ \hline\hline
1 , 3 , 4  &  9193856 \\
1 , 4 , 4  &  307477504 \\
2 , 3 , 3  &  22533736 \\
2 , 3 , 4  &  374111600 \\
2 , 4 , 4  &  12439775296 \\
3 , 4 , 4  &  1810783999024 \\
 & \\
\hline
\end{tabular}
\]
\noindent $\left[\begin{smallmatrix} \frak m\cr \frak n\end{smallmatrix}\right]=\left[\begin{smallmatrix} a&b\cr a&c\end{smallmatrix}\right]:
0\le a<b\le c\le 4$
\[
\begin{tabular}{|c||l||}
\hline
$a,b,c$ & $P_2(a,b;a,c)$\\ \hline\hline
0 , 1 , 1  &  2 \\
0 , 1 , 2  &  4 \\
0 , 1 , 3  &  8 \\
0 , 1 , 4  &  16 \\
0 , 2 , 2  &  8 \\
0 , 2 , 3  &  16 \\
0 , 2 , 4  &  32 \\
\hline
\end{tabular}
\begin{tabular}{|c||l||}
\hline
$a,b,c$ & $P_2(a,b;a,c)$\\ \hline\hline
0 , 3 , 3  &  32 \\
0 , 3 , 4  &  64 \\
0 , 4 , 4  &  128 \\
1 , 2 , 2  &  1684 \\
1 , 2 , 3  &  7112 \\
1 , 2 , 4  &  29200 \\
1 , 3 , 3  &  29728 \\
\hline
\end{tabular}
\begin{tabular}{|c||l|}
\hline
$a,b,c$ & $P_2(a,b;a,c)$\\ \hline\hline
1 , 3 , 4  &  121664 \\
1 , 4 , 4  &  496576 \\
2 , 3 , 3  &  11110312 \\
2 , 3 , 4  &  91632848 \\
2 , 4 , 4  &  754370848 \\
3 , 4 , 4  &  902539626256 \\
 & \\
\hline
\end{tabular}
\]
\noindent $\left[\begin{smallmatrix} \frak m\cr \frak n\end{smallmatrix}\right]=\left[\begin{smallmatrix} a&b\cr c&a\end{smallmatrix}\right]:
0\le a<b\le c\le 4$
\[
\begin{tabular}{|c||l||}
\hline
$a,b,c$ & $P_2(a,b;c,a)$\\ \hline\hline
0 , 1 , 1  &  4 \\
0 , 1 , 2  &  8 \\
0 , 1 , 3  &  16 \\
0 , 1 , 4  &  32 \\
0 , 2 , 2  &  16 \\
0 , 2 , 3  &  32 \\
0 , 2 , 4  &  64 \\
\hline
\end{tabular}
\begin{tabular}{|c||l||}
\hline
$a,b,c$ & $P_2(a,b;c,a)$\\ \hline\hline
0 , 3 , 3  &  64 \\
0 , 3 , 4  &  128 \\
0 , 4 , 4  &  256 \\
1 , 2 , 2  &  1768 \\
1 , 2 , 3  &  7568 \\
1 , 2 , 4  &  31264 \\
1 , 3 , 3  &  32416 \\
\hline
\end{tabular}
\begin{tabular}{|c||l|}
\hline
$a,b,c$ & $P_2(a,b;c,a)$\\ \hline\hline
1 , 3 , 4  &  133952 \\
1 , 4 , 4  &  553600 \\
2 , 3 , 3  &  11205808 \\
2 , 3 , 4  &  92780000 \\
2 , 4 , 4  &  768351424 \\
3 , 4 , 4  &  904335248800 \\
 & \\
\hline
\end{tabular}
\]
\noindent $\left[\begin{smallmatrix} \frak m\cr \frak n\end{smallmatrix}\right]=\left[\begin{smallmatrix} a&a\cr a&b\end{smallmatrix}\right]:
0\le a\le b\le 4$
\[
\begin{tabular}{|c||l||}
\hline
$a,b$ & $P_2(a,a;a,b)$\\ \hline\hline
0 , 0 &  1 \\
0 , 1 &  2 \\
0 , 2 &  4 \\
0 , 3 &  8 \\
0 , 4 &  16 \\
\hline
\end{tabular}
\begin{tabular}{|c||l||}
\hline
$a,b$ & $P_2(a,a;a,b)$\\ \hline\hline
1 , 1 &  82 \\
1 , 2 &  380 \\
1 , 3 &  1624 \\
1 , 4 &  6704 \\
2 , 2 &  151804 \\

\hline
\end{tabular}
\begin{tabular}{|c||l|}
\hline
$a,b$ & $P_2(a,a;a,b)$\\ \hline\hline
2 , 3 &  1303880 \\
2 , 4 &  10791376 \\
3 , 3 &  3300863752 \\
3 , 4 &  54630906416 \\
4 , 4 &  990037617138928 \\
\hline
\end{tabular}
\]
\end{table}


\end{document}